\pdfoutput=1
\documentclass[11pt,a4paper]{scrartcl}

\usepackage[T1]{fontenc}
\usepackage[utf8]{inputenc}
\usepackage{textcomp}
\usepackage{lmodern}
\usepackage{amsmath}
\usepackage{amssymb}
\usepackage{amsthm}
\usepackage{mathtools}
\usepackage{stmaryrd}
\usepackage{bbm}
\usepackage{mathrsfs}
\usepackage{microtype}
\usepackage[shortlabels]{enumitem}
\usepackage[numbers,sort&compress]{natbib}
\usepackage[hidelinks]{hyperref}
\usepackage{bookmark}

\KOMAoptions{DIV=11,BCOR=0mm,parskip=half,headsepline=true}
\addtokomafont{disposition}{\sffamily\bfseries}
\setkomafont{title}{\huge\sffamily\bfseries}
\setkomafont{section}{\Large}
\setkomafont{subsection}{\large}

\mathtoolsset{showonlyrefs}
\numberwithin{equation}{section}

\newcommand*\MSC[1]{\def\mystyleMSC{#1}}
\newcommand*\Keywords[1]{\def\mystyleKeywords{#1}}
\def\mystyleMSC{}
\def\mystyleKeywords{}
\newcommand*\PrintMSCKeywords{%
  \par\smallskip
  \begingroup
  \small
  \noindent\textsc{MSC 2020.} \mystyleMSC\par
  \noindent\textsc{Keywords.} \mystyleKeywords\par
  \endgroup
}
\MSC{82C22; 91D30; 60K35}

\theoremstyle{plain}
\newtheorem{theorem}{Theorem}[section]
\newtheorem{lemma}[theorem]{Lemma}
\newtheorem{proposition}[theorem]{Proposition}
\newtheorem{corollary}[theorem]{Corollary}

\theoremstyle{definition}

\theoremstyle{remark}
\newtheorem{remark}[theorem]{Remark}

\DeclareMathOperator{\supp}{supp}

\newcommand{\bE}{\ensuremath{\mathbb{E}}}

\newcommand{\bN}{\ensuremath{\mathbb{N}}}

\newcommand{\bP}{\ensuremath{\mathbb{P}}}

\newcommand{\bR}{\ensuremath{\mathbb{R}}}

\newcommand{\ind}{\ensuremath{\mathbbm{1}}}

\newcommand{\norm}[1]{\left\Vert \, #1 \, \right\Vert}

\newcommand{\ddx}[1][1]{\ifnum#1=1 \frac{d}{dx} \else \frac{d^{#1}}{dx^{#1}} \fi}
\newcommand{\ddy}[1][1]{\ifnum#1=1 \frac{d}{dy} \else \frac{d^{#1}}{dy^{#1}} \fi}
\newcommand{\ddt}[1][1]{\ifnum#1=1 \frac{d}{dt} \else \frac{d^{#1}}{dt^{#1}} \fi}

\newcommand{\di}{\,\textup{d}}
\newcommand{\bPu}{\bP_\mathbf{u}}

\newcommand{\hP}{\widehat{\bP}}

\newcommand{\bfP}{\mathbf{P}}
\newcommand{\bfx}{\mathbf{x}}

\newcommand{\super}[1]{{(#1)}}

\title{The mean field stubborn voter model}

\author{%
Lisa Hartung\\
\small Institut f\"ur Mathematik, Johannes Gutenberg-Universit\"at Mainz\\
\small \texttt{lhartung@uni-mainz.de}
\and
Christian M\"onch\\
\small Independent researcher\\
\small \texttt{cmoench25@gmail.com}
}
\date{}

\begin{document}

\maketitle

\begin{abstract}
We analyse the effect of agent-dependent heavy-tailed waiting times in the voter model on the complete graph with $N$ vertices. We derive a novel scaling limit and show the existence of a limiting infinite voter model on the slowest updating sites. We further derive the consensus probabilities in the limit model explicitly. In the mean-field setting, the limit is determined by the extreme-value landscape of the waiting times and depends only on the tail index. To obtain these results, we study the coalescing system of random walks that is dual to the limit voter model and prove, among other auxiliary results, that it comes down from infinity.
\end{abstract}
\PrintMSCKeywords


\section{Introduction \& main results}

We investigate a variant of the voter model with agent-dependent update rates, introduced in the (socio-) physics literature by Masuda, Gibert and Redner under the name \emph{heterogeneous voter model} \cite{MasudaGibertRedner10}. As voters with small update rates update slowly, we prefer to refer to this model as the \emph{stubborn voter model}. In the mean-field setting, voters sit on the complete graph with vertices $[N]=\{1,\dots,N\}.$ A voter at site $x\in [N]$ is initially assigned a type (or opinion) $\smash{\eta_0^{(N)}(x)\in\{0,1\}}$ and a clock whose inter-ring times are i.i.d.\ $\smash{\operatorname{Exp}(1/w_x^{(N)})}$, where the mean waiting times $\smash{w_1^{(N)}\geq w_2^{(N)}\geq \dots \geq w_N^{(N)}>0}$ are the order statistics of $N$ i.i.d.\ copies of a positive heavy-tailed random variable $w_1$. We assume that $w_1$ satisfies
\begin{equation}\label{eq:rvtail}
\bP(w_1>t)=t^{-\alpha}L(t),\qquad \alpha\in(0,1),
\end{equation}
where $L$ is slowly varying at infinity. Whenever any clock rings, the corresponding voter copies the opinion of a uniformly chosen neighbour, which gives rise to a Markov process $\smash{(\eta_{t}^{(N)}(x))_{t\geq 0}}$ with right continuous trajectories in $\{0,1\}^{[N]}$. 

We show that a novel scaling limit emerges as the dynamics are dominated by the voters with the highest waiting times. Moreover, the limit depends only on the parameter $\alpha$.

The homogeneous voter model on $\mathbb{Z}^d$ has been extensively studied in the 1970s and 1980s, see \cite{MR0343950,MR0402985,MR0714951,MR0714952} for some important early results and Liggett's monograph \cite{Liggett85} for an overview of the classical theory in the general context of stochastic interacting particle systems. With the advent of Network Science in the last 20 years, interest in models for consensus dynamics has been renewed, particularly for agents that are networked in an inhomogeneous fashion as opposed to the classical translation invariant setting on $\mathbb{Z}^d$. As already mentioned, the finite-population voter model with agent-dependent update rates was introduced by Masuda, Gibert and Redner under the name \emph{heterogeneous voter model} \cite{MasudaGibertRedner10}; related works have studied the influence of `stubborn' or `zealous' agents on opinion dynamics, see for instance \cite{Redner07,yildiz2013binary,Verma14,Waagen14}. Mathematically rigorous results for zealous/stubborn voters were obtained in \cite{DurrettHuo19} and \cite{Roy20}. In the latter models, stubbornness is usually modelled by the presence of a few selected opinions that \emph{never} change, whereas we assume the presence of an inherent heavy-tailed stubbornness distribution that leads to some voters changing their opinion more slowly than others. The resulting model is closer in spirit to trap models and other dynamics in heavy-tailed random environments, except that here the heavy-tailed landscape interacts with an opinion system through the voter-model duality.


We now provide an overview of our main results. We call the sites with the largest weights \emph{slow sites}: a large weight $w_x$ implies a small clock rate $1/w_x$, so these sites take, on average, the longest to update their opinion.
Assume that the common law of $w_1$ satisfies \eqref{eq:rvtail}.
Let $\big(w_i^{(N)}\big)_{i=1}^N$ denote the order statistics of the weights, with $w_1^{(N)}$ being the largest weight. Let $(a_N)_{N\in\bN}$ be any deterministic sequence such that
\begin{equation}\label{eq:aN}
N\bP(w_1>a_N)\longrightarrow 1\qquad \text{as }N\to\infty.
\end{equation}
In the exact Pareto case, one may take $a_N=N^{1/\alpha}$. The slow sites then have weights of order $a_N$. More precisely, let $(\chi_i)_{i\in\bN}$ be i.i.d.\ $\operatorname{Exp}(1)$-distributed and set
$$\xi_i=\Big(\sum_{j\in[i]} \chi_j\Big)^{-\frac{1}{\alpha}}, \quad i\in\bN,$$
then it is well known that
\begin{equation}\label{eq:limweights}
	a_N^{-1}\left(\big(w_i^{(N)}\big)_{i=1}^N, \sum_{i\in[N]}w_i^{(N)} \right) \overset{d}{\longrightarrow} \left((\xi_i)_{i=1}^{\infty}, \sum_{i\in\mathbb{N}}\xi_i\right)\quad \text{ as }N\to\infty,
\end{equation}
see for instance \cite[Chapter 5]{Resni08} and \cite[Section 2.3]{LeadbLindgRootz83}; for the exact Pareto formulation above see also \cite[Theorem 2.33]{Hofst17}. Here and throughout $ \overset{d}{\to}$ denotes convergence in distribution. This convergence of the extreme landscape is the basic mechanism behind our results and is the reason why the limiting process should be viewed as a domain-of-attraction object for heavy-tailed mean-field voter models. Accordingly, the first main result is stated at the level of ranked-environment convergence: once \eqref{eq:limweights} holds, the later arguments no longer use the exact microscopic law.
\begin{theorem}\label{thm:limitvotermodel}
Let $w^{(N)}$ and $\xi$ be as in \eqref{eq:limweights} and let $(\eta_0^{(N)})_{N\in\bN}$ denote a sequence of initial conditions satisfying
\[
\eta_0^{(N)} \overset{d}{\longrightarrow} \eta_0 \text{ as }N\to\infty,
\]
for some random or deterministic $\eta_0\in\{0,1\}^{\bN}$. Then the corresponding Markov processes $(\eta_t^{(N)})_{t\geq 0}$ satisfy
\[
\eta^{(N)}_{a_N(\cdot)} \overset{d}{\longrightarrow} \eta_{(\cdot)} \text{ as }N\to\infty,
\]
where the limit process $\eta$ has initial distribution $\eta_0$ and transition rates governed by $\xi$, i.e.\ conditionally on $\xi$, $(\eta_t)_{t\geq 0}$ is a voter model with update rates $\xi_x,x\in\bN.$
\end{theorem}
\begin{remark}[Remark (Universality, domain of attraction on the complete graph).]
Let $w_1$ satisfy \eqref{eq:rvtail}. Then, after temporal rescaling by $a_N$, the limiting infinite voter model depends only on the tail index $\alpha$ and is independent of the choice of $L$. 
\end{remark}
Note that the existence of the limit process is highly non-trivial, since on the complete graph with vertex set $\bN$, there is no such thing as a `uniformly chosen neighbour' of a vertex. In fact, the assumption of heavy-tailed waiting times is essential for the limit process to exist. Throughout, we work under a fixed environment (quenched), i.e.\ under the distributions
$$\bP=\bP^{\xi}=\mathcal{L}(\eta|\xi)\qquad\text{or}\qquad\bP=\bP^{w}=\mathcal{L}(\eta^{(N)}|w).$$
We next state the consensus probabilities in the limit model.
\begin{theorem}\label{thm:limitvotermodel2}
Let $\eta_0\in\{0,1\}^{\bN}$ be any initial distribution and let
\[\tau=\inf_{t>0}\{\eta_t\equiv 0 \text{ or } \eta_t\equiv 1\}\]
denote the consensus time. Then,
\[
\bE[\tau]<\infty \text{ for almost every realisation of }\xi.
\]
Moreover,
\begin{equation}\label{eq:fixationprobformula}
\bP(\eta_{\tau}\equiv 1)=\frac{\sum_{x\in\bN}\eta_0(x)\xi_x}{\sum_{x\in\bN}\xi_x}.
\end{equation}
\end{theorem}
For completeness, we also state the corresponding formula for the consensus probability in the finite systems.
\begin{theorem}\label{thm:votermodelconsensus}
	Let $\eta_0^{(N)}$ be any initial distribution and let
	\[\tau^{(N)}=\inf_{t>0}\{\eta_t^{(N)}\equiv 0 \text{ or } \eta_t^{(N)}\equiv 1\}\]
	denote the consensus time. 
	Then,
	\[
	\sup_{N} a_N^{-1} \bE[\tau^{(N)}]<\infty \text{ for almost every realisation of }w.
	\]
	Moreover,
	\begin{equation}\label{eq:fixationformula2}
	\bP(\eta_{\tau}^{(N)}\equiv 1)=\frac{\sum_{x\in[N]}\eta^{(N)}_0(x) w_x^{(N)}}{\sum_{x\in[N]}w^{(N)}_x}.
	\end{equation}
\end{theorem}
Note that the consensus time for $\eta^{(N)}$ is trivially finite, because the system has only finitely many states and two of them are absorbing. However, due to the divergence of timescales, the expected consensus times for $\eta^{(N)}$ are not uniformly bounded without rescaling. Below, we prove Theorem~\ref{thm:limitvotermodel2} directly, but an alternative way to derive the asymptotic consensus probabilities is to combine \eqref{eq:fixationformula2} with Theorem~\ref{thm:limitvotermodel} and \eqref{eq:limweights}.  

\begin{remark}[Remark (Extension beyond the complete graph)]
The complete graph is the only setting in which the geometry is fully homogenised and the slow-site landscape is the only source of inhomogeneity. We conjecture that an analogous domain-of-attraction statement should remain valid for broader mean-field graph sequences, at least when the dual random walk mixes on a time scale negligible compared to the extreme waiting-time scale and when its harmonic measure on the set of deepest traps is asymptotically close to uniform.

For asymptotically regular dense graphs, one may therefore expect the same limiting voter model, possibly after a deterministic time change. For more general graph classes, however, the effective motion between slow sites should depend on the graph geometry through a non-trivial trace kernel, so that the correct limit is likely a graph-dependent variant of our infinite voter model rather than a universal one.

Establishing such an extension would require new control of trace processes and meeting probabilities and would lie beyond the scope of the present paper.
\end{remark}

\begin{remark}[Remark (countably many opinions).]
It will be clear from the proofs that Theorems~\ref{thm:limitvotermodel}--\ref{thm:votermodelconsensus} remain valid mutatis mutandis, if opinions are drawn from a countable alphabet instead of $\{0,1\}$. In particular, we may initialise the infinite system in such a way that each voter in $\bN$ carries a distinct opinion. The expected consensus time in this case is still finite and the probability that the whole system ends up with the initial opinion of vertex $x\in\bN$ is $\xi_x/\sum_{z\in\bN}\xi_z$.
\end{remark}

\subsection*{Structure of the paper}
In Section~2, we introduce the system of coalescing random walks that forms the dual of the voter model. Sections~3--5 are devoted to the dual process. In Section~3, we construct the Markov process $X$ that describes the motion of a single walker among slow sites indexed by $\bN$ and derive some of its properties that are of independent interest. In Section~4, we prove an abstract ranked-environment convergence theorem for the walk $Y^{(N)}$ on $[N]$ and extend it to corresponding coalescing systems of finitely many walks. Combined with the extreme-value convergence \eqref{eq:limweights}, this yields the domain-of-attraction statement on the complete graph. In Section~5, we prove Theorem~\ref{thm:comingdown}, which shows that the limit system of coalescing walks comes down from infinity and therefore allows us to treat infinitely many lineages in the voter model simultaneously. Finally, in Section~6, we prove Theorems~\ref{thm:limitvotermodel} and~\ref{thm:limitvotermodel2}.

\section{The dual system: coalescing random walks with heavy tailed rates}
The following duality relation for the voter model is well known, see e.g.\ \cite[Chapter V]{Liggett85}: Given $w^{(N)}$, consider the continuous time random walk $Y^{(N)}$ on $[N]$ with transition rates
$$
q^{(N)}(x,y)=\frac{1}{Nw^{(N)}_x}, \quad x,y\in[N].
$$
Provided with some initial condition $\gamma_0^{(N)}\in \{0,1\}^{[N]}$, which we interpret as the positions of $\sum_{x\in[N]}\gamma^{(N)}_0(x)$ particles, we let each of the particles move independently according to the walk $Y^{(N)}$ and coalesce particles which are located in the same site. This leads to a Markovian system $\big(\gamma_t^{(N)}\big)_{t\geq 0}$ of coalescing random walks that is dual to the voter model in the sense that for all $t\geq 0, A\subset[N]$,
\begin{equation}\label{eq:duality}
\bP\big(\eta^{(N)}_t(x)\equiv 1 \text{ on }A\big| \eta^{(N)}_0\big)=\bP\big(\eta^{(N)}_0(x)\equiv 1 \text{ on }\supp \gamma^{(N)}_t\big| \supp \gamma^{(N)}_0 = A\big),
\end{equation}
where $\supp \gamma_t=\{x:\gamma_t(x)=1\}, t\geq 0.$ Note that the dual dynamics can be reconstructed from the primal dynamics and vice versa. We make use of this fact in Section 6 when proving Theorem~\ref{thm:limitvotermodel} through duality.

\section{Limiting walk on the slow sites}
\subsection*{Construction and properties of the limit process}

The convergence in \eqref{eq:limweights} indicates that a natural choice of metric on the space of environments is the $\ell^1$-distance.
\begin{lemma}
Let $\xi^{(n)},\xi\in\ell_1, n\geq 1,$ be non-negative then $\lim_{n\to\infty} \xi^{(n)}_x=\xi_x$ for all $x\in \bN$ and $$
\lim_{n\to\infty}\sum_{x\in\bN}\xi^{(n)}_x = \sum_{x\in\bN}\xi_x
$$
if and only if $\xi^{(n)}\to\xi$ in $\ell^1$.
\end{lemma}
\begin{proof}
The second assertion obviously implies the first one. For the converse implication, fix $N\in\bN$ and write $Z_{\geq N}(\xi)$ for $\sum_{x\geq N}\xi_x$, then
\[
|\xi-\xi^{(n)}|_{\ell^1}\leq \sum_{x<N}|\xi_{x}^{(n)}-\xi_x| + Z_{\geq N}(\xi)+ Z_{\geq N}(\xi^{(n)}).
\]
Fix $\epsilon>0$, then $N$ so large that $Z_{\geq N}(\xi)\leq \epsilon$. Fix $L$ so large that $|Z_{\geq N}(\xi)-Z_{\geq N}(\xi^{(n)})|\leq \epsilon$ whenever $n\ge L$. We obtain for $n\ge L$
\[
|\xi-\xi^{(n)}|_{\ell^1}\leq \sum_{x<N}|\xi_{x}^{(n)}-\xi_x|+3\epsilon,
\]
and taking the limit $n\to\infty$ yields $|\xi-\xi^{(n)}|_{\ell^1}<3\epsilon$ due to coordinate-wise convergence, which concludes the argument since $\epsilon$ is chosen arbitrarily.
\end{proof}
Let $\bfP$ denote the distribution of the limiting environment $\xi$. Note that the strong law of large numbers yields that $\bfP$-almost surely
$$
i/2 < \sum_{j\in[i]}\eta_i <2i, \text{ for all but finitely many }i\in\bN,
$$
and consequently we have
\begin{equation}\label{eq:xipprox}
	\xi_i \asymp i^{-1/\alpha}\quad \bfP\text{-almost surely}.
\end{equation}
Here, $f(i)\asymp g(i)$ means that $f(i)/g(i)$ is bounded away from $0$ and $\infty$ uniformly in $i\in\bN$. In particular, by \eqref{eq:xipprox}, the right-hand side of \eqref{eq:limweights} is a well-defined random variable on $\ell_1\times (0,\infty)$. We have that
$$
	\xi([N])
=\sum_{x\in[N]}\xi_x\overset{N\to\infty}{\longrightarrow}\xi(\mathbb{N})=\sum_{x\in\mathbb{N}}\xi_x<\infty,\quad \bfP\text{-almost surely},$$
and that 
$$
\label{eq:xilambdasum}
Z(\lambda)=\sum_{x\in\bN}\frac{\lambda \xi_x}{1+\lambda\xi_x}\in\bR,\quad \text{ exists }\bfP\text{-almost surely for any }\lambda\in\bR\setminus\{-\xi_{x}^{-1},x\in\bN\}.
$$
Eventually, our goal is to construct a limiting coalescing system that captures the behaviour of the coalescing random walk amongst the slow sites, thereby extending \eqref{eq:limweights} from the environments to the processes. Our first goal is thus to construct a random walk on $\bN$ whose transitions are governed by the limit environment $\xi=(\xi_i)_{i\in\bN}$ described prior to \eqref{eq:limweights}. As the holding times outside the slow sites tend to zero under the rescaling, this construction must be done carefully.
We start the construction with finitely many sites. Let $X^{(N)}=(X_t^{(N)})_{t\geq0}$ be a random walk on $[N]$ with uniform jump distribution and waiting times governed by $\xi$. We assume that the starting point $X_0^{(N)}\in [N]$ is chosen according to some distribution $\nu_N$ with $\supp(\nu_N)\subset[N]$. For technical convenience, we consider $X^{(N)}$ as a Markov process on $\bN$ via the canonical embedding of $[N]$ into $\bN$, which leads to the transition rates
$$
r^{(N)}(x,y)=\begin{cases}
	 \frac{1}{N\xi_x}, & \text{ if }x,y\in[N],\\
	 0, & \text{ otherwise.}
\end{cases}
$$
\begin{remark}[Remark.]
Note that at every jump time $X^{(N)}$ may remain at its current position with probability $1/N$, i.e.\ it is a lazy version of the random walk on $[N]$ with holding times governed by $\xi$. Defining $X^{(N)}$ in this manner leads to slightly simpler calculations. It is easy to see that the limit process (see  Theorem~\ref{thm:OSTRW} below) is the same for both the lazy version $X^{(N)}$ and the non-lazy version.
\end{remark}
We mostly work under the quenched measure $\bP=\bP^\xi$ conditional on the environment $\xi$. We use the shorthand $\bP_x(\,\cdot\,)$ for $\bP(\,\cdot \,|X^{[N]}_0=x),$ where it is tacitly assumed that $N\geq x$. Similarly we define $\bP_{\nu^N}=\sum_{x\in[N]}\nu^N(x)\bP_x$, and write $\bP_\nu$ instead of $\bP_{\nu^N}$ when there is no ambiguity regarding the initial distribution. We frequently initiate $X^{(N)}$ in a uniformly chosen site and denote the corresponding distributions by $\bPu=N^{-1}\sum_{x\in[N]}\bP_x.$

Note that $X^{(N)}$ only depends on the initial segment of $\xi$. If we keep one realisation of $\xi$ fixed and choose $M>N$, then $X^{(N)}$ can be viewed as a projection of $X^{(M)}$ up to a time change. To make this precise, let $W=(W_t)_{t\geq 0}$ be any Markov process on $\bN$ under $\bP$ with stationary distribution $\nu$. Let $A\subset \bN$ be such that $\nu(A)>0$, thereby ensuring that, $\bP$-almost surely,
$$
W^{-1}(A)=\{t: W_t\in A \}
$$
is unbounded and $T_A=\inf\{t>0:W_t\in A\}$ is finite. Let now $\psi^{W}_{A}:[0,\infty)\to W^{-1}(A)$ be given by
$$ \psi^W_{A}(t)=\int_{0}^{t} \ind_{\{ W_t\in A \}}\,\textup{d}t.$$ 
Since $\psi_A^{W}$ is non-decreasing, it admits an increasing right-continuous inverse $\phi_A^W$. We call the time-changed process $\textup{tr}_A(W)$ defined via
$$
\textup{tr}_A(W)_t=W_{\phi^W_A(t)},\quad t\geq 0,
$$
the \emph{trace} of $W$ on $A.$
\begin{lemma}\label{lem:project}
Let $M>N\in\bN$. For any $x\in[N]$, the processes $\textup{tr}_{[N]}\big(X^{(M)}\big)$ and $X^{(N)}$ have the same distribution under $\bP_x$. The same assertion remains true, if both $X^{(M)}$ and $X^{(N)}$ are started in a uniformly chosen site in $[M]$ and $[N]$, respectively.
\end{lemma}
\begin{proof}
For any given realisation of $\xi$, this follows immediately from the definition of $\textup{tr}_{[N]}\big(X^{(M)}\big)$ and the uniformity of the jump distributions of $X^{(N)}$ and $X^{(M)}$.
\end{proof}
The next lemma establishes that the time spent outside the slow sites is asymptotically negligible.
\begin{lemma}\label{lem:notimeoutside}
For $M> N\geq 1$ let $\phi_{M,N}:[0,\infty)\to[0,\infty)$ be the associated time change of the trace process $\textup{tr}_{[N]}\big(X^{(M)}\big)$ of $X^{(M)}$ started in a fixed initial distribution $\nu^M$ on $[M]$. For any sequence $(\epsilon_N)_{N\in\bN}$ in $(0,1)$, with $$
\lim_{N\to\infty}\epsilon_N N^{\frac1\alpha-1}=\infty,
$$
it holds that
$$
\lim_{N\to\infty}\sup_{M > N}\bP_{\nu^M}(|\phi_{M,N}(t)-t|\geq \epsilon_N)=0, \quad t>0.
$$
\end{lemma}
Throughout the proof of Lemma \ref{lem:notimeoutside} and the remainder of this section, we frequently decompose the path of $X^{(M)}$ into excursions. To prepare the arguments, we establish a few auxiliary results pertaining to the excursion structure of $X^{(M)}$ first. For a given proper non-empty subset $A$ of $[M]$, denote by $T_0^A$ the first hitting time of $A$ by $X^{(M)}$ and by $\mathcal{E}^{A}_1,\mathcal{E}^{A}_2,\dots$ the successive excursion intervals of $X^{(M)}$ from $A$ after time $T^A_0$. Note that $T^A_j<\infty$ for all $j\in\bN_0$ $\bP$-almost surely and that $T^A_0=0$ precisely if $X_0\in A.$ Clearly, the set of times $[T^A_0,\infty)\setminus \bigcup_{n=1}^{\infty}\mathcal{E}^A_n$ can also be identified with a sequence $\mathcal{I}^A_1,\mathcal{I}^A_2,\dots$ of non-adjacent $\bP$-almost surely finite intervals. Let $E^A$ have the distribution of $|\mathcal{E}^A_1|$ and let $I^A$ have the distribution of $|\mathcal{I}^A_1|$. By stationarity of the transition rates, uniformity of the jump distribution and the Markov property, we have that
$$
({E}^A_n)_{n\in\bN}:=(|\mathcal{E}^A_n|)_{n\in\bN} \text{ and } ({I}^A_n)_{n\in\bN}:=(|\mathcal{I}^A_n|)_{n\in\bN}
$$
are two independent sequences of i.i.d.\ copies of $E^A$ and $I^A$. We refer to $E^A$ and $I^A$ as the generic \emph{excursion length} and the generic \emph{incursion length} of $A$, respectively.  Observe also, that $T_0\overset{d}{=}E^A$ if and only if $X_0$ is distributed uniformly on $[M]\setminus A$. We further define $H^{A}$ to be $\operatorname{Exp}(\xi_U^{-1})$-distributed conditionally on $U$, where $U$ is uniform on $A$. Hence, $H^{A}$ is a \emph{typical holding time} of $X^{(M)}$ conditioned to stay in $A$. The next two lemmas collect some properties of holding times and excursion/incursion lengths.
\begin{lemma}\label{lem:LD_holding}
Let $A\subset[M]$ be such that $\max_{x\in A}\xi_x\leq 1/2$ and let $H^A_1,H^A_2,\dots$ be i.i.d.\ typical $A$-holding times. Set
\[
\bar\xi(A)=\frac{\sum_{x\in A}\xi_x}{|A|}.
\]
For any $l\in\bN$ and $\lambda>0$, we have
\[
\bP\left(\sum_{j=1}^l H^A_i>l\lambda\right)\leq \textup{e}^{-l(\lambda-2\bar{\xi}(A)+O(\bar{\xi}(A)^2))}.
\]
\end{lemma}
\begin{proof}
Markov's inequality together with the fact that $\bE \textup{e}^{H^A_1}=|A|^{-1}\sum_{x\in A}\frac{1}{1-\xi_x}$ yields
\[
\bP\left(\sum_{j=1}^l H^A_i>l\lambda\right)=\bP\left(\textup{e}^{\sum_{j=1}^l H^A_i}>\textup{e}^{l\lambda}\right)\leq \textup{e}^{-l\lambda + l \log(|A|^{-1}\sum_{x\in A}\frac{1}{1-\xi_x})}.
\]
Now the assertion follows by a Taylor expansion of $\log(\cdot)$ around $1$, since
\begin{equation}\label{eq:xibound}
\frac{1}{|A|}\sum_{x\in A}\frac{1}{1-\xi_x}=1+\frac{1}{|A|}\sum_{x\in A}\frac{\xi_x}{1-\xi_x}\leq 1+2\bar{\xi}(A)\leq2,
\end{equation}
by assumption on $(\xi_x)_{x\in A}$.
\end{proof}

\begin{lemma}
Fix an i.i.d.\ sequence of typical $A$-holding times $(H^A_i)_{i\in\bN}$ and let $G$ denote a $\operatorname{Geo}(1-|A|/M)$ random variable independent of $(H^A_i)_{i\in\bN}$. Then,
\begin{equation}\label{eq:increp}
I^A \overset{\textup{d}}{=} \sum_{i=1}^{G}H^{A}_i,
\end{equation}
and furthermore
\begin{equation}\label{eq:Ulaplace}
	\bE\Big[\textup{e}^{-\lambda I^A} \Big] = \frac{(M-|A|) \, \frac{1}{|A|}\sum_{x\in A} \frac{1}{1+\lambda \xi_x}}{ M-\sum_{x\in A}\frac{1 
		}{1+\lambda \xi_x}},\quad \lambda\in\bR\setminus\{-\xi_x^{-1}, x\in A\},
\end{equation}
as well as
\begin{equation}\label{eq:UVexpectation}
	\bE I^A = \frac{M}{M-|A|}\,\frac{1}{|A|} \sum_{x\in A}\xi_x.
\end{equation}
In particular, for $A=[N]$ with $N<M$, it follows
that $\bfP$-almost surely for $E^{[N]}=I^{[M]\setminus[N]}$
\begin{equation}\label{eq:V1decay}
	\sup_{M>N}\bE \big[E^{[N]}\big] \in O\big( N^{-\frac1\alpha} \big).
\end{equation}
\end{lemma}
\begin{proof}
The representation \eqref{eq:increp} follows from two simple observations. Firstly, upon entering $A$, $X^{(M)}$ sits in a uniformly distributed site and has a site-independent probability of $1-\frac{|A|}{M}$ at each jump to leave $A$ and each jump within $A$ leads with probability $|A|^{-1}$ to any given site $x$. Secondly, at each site $x$, the holding time is an independent $\operatorname{Exp}(1/\xi_x)$-random variable. \eqref{eq:Ulaplace} can then be obtained from \eqref{eq:increp} by elementary calculations using 
\begin{equation}\label{eq:moment}
\bE \textup{e}^{-\lambda H^A_1}=|A|^{-1}\sum_{x\in A}\frac{1}{1+\lambda\xi_x},
\end{equation}
 and the independece of $G$ and $(H^A_i)_{i\in\bN}$. \eqref{eq:UVexpectation} follows directly from \eqref{eq:Ulaplace}.  \eqref{eq:V1decay} is a straightforward consequence of \eqref{eq:UVexpectation} and \eqref{eq:xipprox}.
\end{proof}
Next, we prove Lemma \ref{lem:notimeoutside}.
\begin{proof}[Proof of Lemma \ref{lem:notimeoutside}]
We consider the excursions of $X^{(M)}$ from $[N]$. Set $T_j=T_j^{[N]}$, $j\in\bN_0$. Let, for fixed $t>0$, $$\ell(t)=\max\{n\geq 0 : T_n \leq \phi_{M,N}(t) \}.$$ Note that for any $t\geq 0$, $$\phi_{M,N}(t)-t = T_0+\sum_{j=1}^{\ell(t)}E_j \geq 0,$$
with $E_j, j\in\bN$ i.i.d.\ copies of $E^{[N]}$. In particular, for any fixed $\epsilon>0$
\begin{equation}\label{eq:split}
\bP_{\nu}(\phi_{M,N}(t)-t>2\epsilon)\leq \bP_{\nu}(T_0\geq \epsilon) + \bP\Big( \sum_{j=1}^{\ell(t)}E_j >\epsilon \Big).
\end{equation}
For the second term on the right hand side of \eqref{eq:split}, we first note that it does not depend on the initial condition \eqref{eq:Ulaplace}. Moreover,
\begin{equation}\label{eq:decomp}
\bP\Big( \sum_{j=1}^{\ell(t)}E_j >\epsilon \Big)\leq \bP\Big( \sum_{j=1}^{L}E_j >\epsilon \Big)+\bP(\ell(t)>L),\quad L\in\bN.
\end{equation}
We denote the incursion lengths inside $[N]$ by $I_1,I_2,\dots$. Then, by Chebyshev's inequality 

$$
\bP(\ell(t)>L)\leq \bP\Big(\sum_{j=1}^L I_j<t\Big)\leq \textup{e}^{\lambda t}\bE\Big[\textup{e}^{-\lambda I_1} \Big]^L\leq \textup{e}^{-q_N(\lambda)L+\lambda t},\quad L\in\bN, \lambda>0,
$$
where we used that by \eqref{eq:Ulaplace}, 
\begin{equation}\label{eq:weakbd}
\bE\Big[\textup{e}^{-\lambda I_1} \Big] = \frac{(M-N)(1-q_N(\lambda))}{M-N(1-q_N(\lambda))}\leq 1-{q_N(\lambda)}, \quad \lambda>0,
\end{equation}
with $q_N(\lambda)\equiv1-N^{-1}\sum_{x\in\bN}\frac{1 
}{1+\lambda \xi_x}$. 
Fixing $K_0\in\bN$, we obtain for any choice of $L=L(N,M)$ satisfying $L>K_0 N$
\begin{equation}\label{eq:LDBD}
\lim_{N\to\infty}\sup_{M>N}\bP(\ell(t)>L)\leq \textup{e}^{-K_0 Z(\lambda)+\lambda t}, \quad \lambda >0,
\end{equation}
since $N q_N(\lambda)=\sum_{y\in[N]}\frac{\lambda \xi_y}{1+\lambda \xi_y}\to Z(\lambda)$ as $N\to\infty.$ On the other hand, for $L=2K_0N$, we bound the second term in \eqref{eq:decomp} from above by
\begin{equation}\label{eq:domterm}
\bP\Big( \sum_{j=1}^{L}E_j >\epsilon \Big)\leq \frac{2K_0N \bE E_1}{\epsilon}\leq \frac{C K_0}{\epsilon}N^{1-1/\alpha}
\end{equation}
for some $C$ depending only on $\xi$, by Markov's inequality and \eqref{eq:V1decay}. We may thus conclude from \eqref{eq:decomp} and \eqref{eq:LDBD} that
\begin{equation}\label{eq:final1}
\lim_{N\to\infty}\sup_{M>N}\bP\Big( \sum_{j=1}^{\ell(t)}E_j >\epsilon \Big)=0,
\end{equation}
since $K_0$ can be made arbitrarily large. It remains to show that $\bP_\nu(T_0>\epsilon)$ vanishes. If $X_0\in[N]$, then $T_0=0$. Otherwise, the first jump of $X$, say at time $\sigma_0$, either sends $X$ to $[N]$ in which case $\sigma_0=T_0$ or to a uniformly distributed point in $[M]\setminus[N]$ in which case the remaining time to enter $[N]$ has precisely the distribution of $E_1$. We conclude
$$
\bP_\nu(T_0>\epsilon)\leq \sum_{x\in[M]\setminus[N]}\bP_x(\sigma_0>\epsilon/2)\nu(x) +\bP(E_1>\epsilon/2)
$$
and thus 
\begin{equation}\label{eq:final2}
	\lim_{N\to\infty}\sup_{M>N}\bP_{\nu}(T_0>\epsilon)\leq \lim_{N\to\infty}\sup_{M>N} \big(\bP(\tau_{N+1}>\epsilon/2)+\bP(E_1>\epsilon/2)\big) = 0,
\end{equation}
where we have used that $\sigma_0$ is dominated by an $\operatorname{Exp}(1/\xi_{N+1})$ random variable $\tau_{N+1}$ uniformly in the starting point $x\in[M]\setminus[N]$ and that $\sup_{M>N}\bP(E_1>\epsilon/2)\in O(\epsilon^{-1}N^{-1/\alpha})$ by \eqref{eq:V1decay} and Markov's inequality. Combining \eqref{eq:final1} and \eqref{eq:final2} according to \eqref{eq:split} concludes the proof for fixed $\epsilon>0$. The assertion for decaying $\epsilon_N$ follows from the observation that \eqref{eq:LDBD} holds independently of the choice of $\epsilon$ and that in the remaining estimates \eqref{eq:domterm} and \eqref{eq:final2} the term on the right hand side of \eqref{eq:domterm} dominates. 
\end{proof}
We next describe the limiting process associated with the sequence $(X^{N})_{N\in\bN}$. For $N\in\bN$, let $P^{(N)}=(P^{(N)}_t)_{t\geq0}$ denote the semigroup of $X^{(N)}$ given $\xi.$
\begin{proposition}\label{prop:Qexistence}
Let $(\nu^M)_{M\in\bN}$ be distributions with $\supp\nu^M\subset[M]$ for each $M\in\bN$ and $\nu^M([N])\to0$ for every $N\in\bN$. Then there exists a family of probability laws $(Q_t)_{t>0}$ on $\bN$, independent of $(\nu^M)_{M\in\bN}$, such that
$$
	Q_t = \lim_{M\to\infty} \nu^M P^{(M)}_t,\; t>0,\quad \bfP\text{-almost surely.}
$$
\end{proposition}
Before we give the proof, we formulate an auxiliary result that is also useful later on.
\begin{lemma}\label{lem:Tailbound}
For any $M>N>2$ and $t>0$, we have that
\begin{equation*}\begin{aligned}
 \bPu(X^{(M)}_t\in \bN\setminus[N] & |X_0^{(M)}\in[M]\setminus[N])\\
& \leq \frac{2M-N}{M-N}\left( \textup{e}^{-s_Nt}+\frac{N}{s_N}\left(\frac{t}{\xi([N])}+ C\right)\right),
\end{aligned}
\end{equation*}
where $C$ is a universal constant and
$$s_N=s_N(\xi)=\frac{1}{2}\left(1/\xi_{N+1} \wedge N/2\xi(\bN\setminus[N])\right).$$

\end{lemma}
\begin{proof}
Let $I=I^{[N]}$ and $E=E^{[N]}$  be the generic incursion and excursion lengths of $X^{(M)}$ from $[N]$. For any $M>N$, we set
$$
Z_t^{(M)}=\ind_{\{ X^{(M)}\in [M]\setminus[N] \}}, \quad t\geq 0.
$$
Conditionally on $X_0^{(M)}\in[M]\setminus [N]$, the starting point $X_0^{(M)}$ is uniformly distributed on $[M]\setminus [N]$ and thus the conditional distribution of $(Z_t^{(M)})_{t\geq 0}$ is that of a simple two state renewal process, with transition function $p^M(r)=\bPu(Z_r^{(M)}=1|X_0^{(M)}\in[M]\setminus [N])$ given by
\begin{equation}\label{renewal}
p^M(r)=\bar{F}_{E}(r)+\int_{0}^r \bar{F}_{E}(r-s) \di \ell^{M}(s), r\geq 0.
\end{equation}
Here, $\bar{F}_{E}(r)=\bP(E>r)$ and $\ell^{M}$ is the renewal measure associated with the generic cycle length $E+I$. Let us first provide bounds on $\bar F_E$. Let $G$ denote the number of stays inside $[M]\setminus [N]$ that $X^{(M)}$ accumulates during a generic excursion $E$, and let $H_1,H_2,\dots$ denote independent generic $[M]\setminus [N]$-holding times. 
We have, by Chebyshev's inequality,
\begin{align*}
	\bP(E>t) &= \bP\left(\sum_{j=1}^{G} H_j>t\right) = \sum_{j=1}^{\infty}\frac{N}{M}(1-N/M)^{j-1}\bP\left(\sum_{i=1}^j H_i>t \right) \\ 
	&\leq \sum_{j=1}^{\infty}\frac{N}{M}(1-N/M)^{j-1}\textup{e}^{-st}(\bE\textup{e}^{s H_1})^j.
\end{align*} 
Using \eqref{eq:moment} and taking $s$ equal to   $$s_N=s_N(\xi)=\frac{1}{2}\left(1/\xi_{N+1} \wedge N/2\xi(\bN\setminus[N])\right),$$
which is positive $\bfP$-almost surely for all $N>2$, we get
\begin{align*}
	\bP(E>t) &\leq \frac{N}{M-N} \textup{e}^{-s_Nt}\sum_{j=1}^{\infty}\left(\frac{1}{M}\sum_{x=N+1}^M \frac{1}{1-s_N\xi_x} \right)^j\\
	&\leq \frac{N}{M-N} \textup{e}^{-s_Nt}\sum_{j=1}^{\infty}\left(\frac{M-N+2s_N\xi(\bN\setminus[N])}{M} \right)^j\\
	& \leq \textup{e}^{-s_Nt} \frac{N}{M-N} \frac{2M-N}{N}=\textup{e}^{-s_Nt}\frac{2M-N}{M-N}.
\end{align*} 
The second inequality uses a bound similar to \eqref{eq:xibound}. The third inequality is also true by the choice $s_N$.  Inserting this estimate into \eqref{renewal} yields
\begin{equation}\label{eq:Rbd1.1}
	p^M(t) \leq \frac{2M-N}{M-N}\left( \textup{e}^{-s_Nt}+\int_{0}^t \textup{e}^{-s_N(t-s)} \di \ell^{M}(s)\right).
\end{equation}
Note that $\ell^M$ is dominated by the renewal measure $r^N$ associated with the typical holding time $H=H^{[N]}$ inside $[N]$, because each cycle contains at least one stay inside $[N]$. We obtain 
$$
\int_{0}^t \textup{e}^{-s_N(t-s)} \di \ell^{M}(s)\leq \frac{r^N(t)}{s_N}\int_{0}^t s_N \textup{e}^{-s_N(t-s)} \frac{1}{r^N(t)}\di r^{N}(s)\leq \frac{r^N(t)}{s_N},
$$
where the last inequality follows from the fact that the integral can be viewed as the distribution function at $t$ of the sum of two independent random variables. We conclude that
\begin{equation}\label{eq:Rbd1.2}
	p^M(t)\leq \frac{2M-N}{M-N}\textup{e}^{-s_Nt}+\frac{2M-N}{M-N}\frac{r^N(t)}{s_N},
\end{equation}
and the renewal function $r^N(t)$ satisfies the bound
$$
r^N(t)\leq \frac{t}{\bE H}+ C \frac{\bE H^2}{(\bE H)^2}= \frac{tN}{\xi([N])}+ C\frac{2N\sum_{x=1}^N \xi_x^2}{\xi([N])^2},
$$
where $C$ is some universal constant, see e.g.\ \cite{Daley76}. Combining this bound with \eqref{eq:Rbd1.2} concludes the proof.
\end{proof}

\begin{proof}[Proof of Proposition~\ref{prop:Qexistence}]
We set
$$R_t^{M}(y) = \frac{1}{M}\sum_{x\in[M]} \delta_x P^{(M)}_t(y)=\bPu(X^{(M)}_t=y),  \quad y\in[M],\, M\in\bN,$$
and divide the proof into three parts: first we establish tightness of $(R_t^{(M)})_{M\in\bN}$, then we show the existence of a unique accumulation point, and finally we argue that the uniform distributions in the definition of $R^{(M)}_t$ may be replaced by any sequence $(\nu^M)_{M\in\bN}$ that puts vanishing mass on the slow sites.\\
\noindent\textit{Part (i) -- tightness:} Let $\epsilon>0$ be given and fix $N\in \bN$. Let now $M_0>N/\epsilon$. We obtain
\begin{equation}\label{eq:straff1}
R_t^{M}(\bN\setminus[N])\leq \epsilon + \frac{M-N}{M}\bPu(X^{(M)}_t\in \bN\setminus[N]|X_0^{(M)}\in[M]\setminus[N]), \quad M>M_0, t>0.
\end{equation}
Thus, Lemma~\ref{lem:Tailbound} implies that for some constant $C$ depending only on $\xi$, and all sufficiently large $M$,
\begin{equation}\label{eq:ConeandCtwo}
R_t^{M}(\bN\setminus[N]) \leq  \epsilon + C t N^{1-1/\alpha},
\end{equation}
since $1/s_N=O(N^{-\frac1\alpha})$ $\bP$-alomst surely. It follows that $(R_t^{M})_{M\in\bN}$ is indeed tight.\\ 
\noindent\textit{Part (ii) -- uniqueness:} We use a renewal theoretic argument to obtain the uniqueness of the accumulation point. Fix $y\in \bN$ and let
$$
W_s^{(M)}=\ind_{\{ X^{(M)}_s=y\}}, \quad s\geq 0.
$$
Since, for each $y\in\bN$, $\bP_{\mathbf{u}}(X^{(M)}_s=y)=\bP(W_s^{(M)}=1)$, it suffices to show that the sequence $(W^{(M)})_{M\geq 1}$ of renewal processes converges in distribution, since this identifies the limits of $(R_t^{(M)}(y))_{M\in\bN}, y\in\bN$. Furthermore, on the event $\{X_0^{(M)}\neq y\}$, each $W^{(M)}$ is entirely characterised by the sums
\[
\sum_{j=1}^{n} E_j\ind_{j \text{ odd}}+I_j\ind_{j \text{ even}},\quad n\in\bN, 
\]
where $E_j,j\in\bN$ are i.i.d.\ copies of $E(M)=I^{[M]\setminus\{y\}}$ and $I_j,j\in\bN$ are i.i.d.\ copies of $I(M):=I^{\{y\}}$, respectively. Since $\lim_{M\to\infty}\bP_{\mathbf u}(X_0^{(M)}\neq y)=1$, it is thus easy to see that convergence in distribution of $E(M), I(M)$ as $M\to\infty$ implies the convergence of $(W^{(M)})_{M\in\bN}$. By \eqref{eq:increp}, $I(M)$ is the independent sum of $G\sim\operatorname{Geo}(1-1/M)$ copies of $\operatorname{Exp}(1/\xi_y)$, thus its convergence to $\operatorname{Exp}(1/\xi_y)$ is immediate. To obtain the limit of $E(M)$, we use the corresponding characteristic function
\[
\psi_M(\lambda)=\bE\left[\textup{e}^{i \lambda E(M)}\right]=\frac{\frac{1}{M-1}\sum_{x\in [M]\setminus\{y\}} \frac{1}{1-i\lambda \xi_x}}{ M-\sum_{x\in [M]\setminus\{y\}}\frac{1}{1-i\lambda \xi_x}},\quad \lambda\in\bR,
\]
which is obtained in the same manner as \eqref{eq:Ulaplace}. The pointwise limit of $(\psi_M)_{M\in\bN}$ is simply
\begin{equation}\label{eq:exccharf}
\psi(\lambda)=\frac{1}{1-\sum_{x\in \bN\setminus\{y\}} \frac{i\lambda \xi_x}{1-i\lambda \xi_x}},\quad \lambda\in\bR,
\end{equation}
which is a well-defined point on the complex unit circle since $|1-i\lambda \xi_x|\geq 1$ and thus $\sum_{x\in \bN\setminus\{y\}} \left|\frac{i\lambda \xi_x}{1-i\lambda \xi_x}\right|\leq \lambda \xi(\mathbb{N}) <\infty$, for any $\lambda\in\bR.$ L\'evy's continuity theorem now yields that $\psi$ is the characteristic function of a probability distribution if it is continuous in $0$, which follows from
\[
\left|\frac{1}{1-\sum_{x\in \bN\setminus\{y\}} \frac{i\lambda \xi_x}{1-i\lambda \xi_x}}-1\right|\leq \frac{\left|\sum_{x\in \bN\setminus\{y\}} \frac{i\lambda \xi_x}{1-i\lambda \xi_x}\right|}{\left|1-\sum_{x\in \bN\setminus\{y\}} \frac{i\lambda \xi_x}{1-i\lambda \xi_x}\right|}\leq \lambda \xi(\bN)|\psi(\lambda)|=\lambda \xi(\bN). 
\]
We conclude, that $E((M))_{M\in\bN}$ and $(I(M))_{M\in\bN}$ both have well-defined distributional limits, whence $(W^{(M)})_{M\in\bN}$ converges in distribution to a well-defined limiting renewal process.
\\
\noindent\textit{Part (iii) -- independence of initial distribution:} Let $\tau$ denote the first time that ${X}^{(M)}$ jumps, then for any $y\in [M]$
$$
\bP_{\nu^M}({X}^{(M)}_{\tau+s}=y)=\bP_{\mathbf{u}}({X}^{(M)}_{s}=y)=R_s^M(y), \quad s>0,
$$
and, using the Markov property,
\begin{align}\label{eq:indep1}\bP_{\nu^M}({X}^{(M)}_{t}=y, t>\tau) & =\bP_{\nu^M}({X}^{(M)}_{\tau+(t-\tau)}=y,t>\tau)\\
	& =\sum_{x\in[M]}\nu^M(x)\int_{0}^t R_{t-s}^M(y) \frac{\textup{e}^{-s/\xi_x}}{\xi_x}\di s.
\end{align}
By the assumption that $\nu^M([N])$ vanishes for any fixed $N$, we obtain that $$\bP(\tau>t)\leq \nu^M([N])+\int_t^\infty \frac{\textup{e}^{-s/\xi_N}}{\xi_N},$$
which can be made arbitrarily small by choosing $N$ large and then increasing $M$. Hence $\lim_{M\to\infty}\bP_{\nu^M}(\tau<t)=1$. On the other hand, we have that
\begin{align*}\sum_{x\in[M]\setminus[N]}\nu^M(x)\int_{0}^t & R_{t-s}^M(y) \frac{\textup{e}^{-s/\xi_x}}{\xi_x}\di s    \leq \sum_{x\in[M]}\nu^M(x)\int_{0}^t R_{t-s}^M(y) \frac{\textup{e}^{-s/\xi_x}}{\xi_x}\di s\\
	&  \leq \sum_{x\in[M]\setminus[N]}\nu^M(x)\int_{0}^t R_{t-s}^M(y) \frac{\textup{e}^{-s/\xi_x}}{\xi_x}\di s+ \nu^M([N]),
\end{align*}
 and the weak convergence of $\operatorname{Exp}(\lambda)$ to $\delta_0$ as $\lambda\to\infty$ implies that both the first and the last term can be brought arbitrarily close to $R_t^M(y)$ by choosing first $N$ and then $M$ large. It thus follows from \eqref{eq:indep1} that
 $$\lim_{M\to\infty}\bP_{\nu^M}(X_t = y)=\lim_{M\to\infty}R_t^M(y)=Q_t(y),\quad y\in \bN, t>0.$$
\end{proof}
The laws $(Q_t)_{t>0}$ describe the limiting distribution of the processes $X^{(N)}, N\in\bN$ upon (re-)entering the slow sites. To define the dynamics of the limiting process, let $\tau_x, x\in\bN$, be $\operatorname{Exp}(1/\xi_x)$ distributed and set, for $f\in\mathcal{B}(\bN)=\{g:\bN \to \bR: g \text{ bounded }\}$,
\begin{align}\label{eq:P_t-decomposition}
	P_tf(x) = f(x)\bP(\tau_x>t) + \bE(Q_{t-\tau_x}(f)\ind_{\{\tau_x<t\}}),\quad x\in\bN, t\geq 0.
\end{align}
\begin{theorem}\label{thm:OSTRW}
The following assertions hold $\bfP$-almost surely:
\begin{enumerate}[(a)]
	\item We have $\lim_{N\to\infty}P_t^{(N)}=P_t$ for every $t>0.$
	\item $P=(P_t)_{t\geq 0}$ is a Markov semigroup inducing a process $X=(X_t)_{t> 0}$ on $\bN$ with waiting times governed by $\xi$.
	\item $(Q_t)_{t\geq0}$ is an entrance law of $X$, i.e.\ satisfies $Q_t P_s = Q_{t+s}$, for all $s,t>0.$
	\item The trace process $\textup{tr}_{[N]}(X)$ coincides with $X^{(N)}$ in distribution.
\end{enumerate}	
\end{theorem}
\begin{proof}
We first prove \textit{(a)}. Fix any bounded function $f$ on $\bN$ and any $t>0$. We have $$P_t^{(N)}f(x)=f(x)\bP(\tau_x>t)+ \bE(R^{N}_{t-\tau_x}(f)\ind_{\{\tau_x<t\}}),\quad x\in[N],$$
and it thus follows from Proposition~\ref{prop:Qexistence} and dominated convergence that \begin{equation}\label{eq:Plimit}\lim_{N\to\infty}P_t^{(N)}f(x)=P_tf(x),
\end{equation}
for any $x\in\bN.$  To show that the identity in \textit{(c)} holds, we observe that 
$$
R^N_tP^{(N)}_s=R^{N}_{t+s}, \quad t,s\geq 0.
$$
\eqref{eq:Plimit} together with Proposition~\ref{prop:Qexistence} allows us to take weak limits and establish the identity
\begin{equation}\label{eq:entrance}Q_tP_s=Q_{t+s}, \quad t,s\geq 0.
\end{equation}
Next, we show the semigroup property needed for \textit{(b)}. \textit{(c)} implies, for any $x\in\bN$,
\begin{align*}
 P_{t+s}f(x) =\, & \bP(\tau_x>t+s)+ \bE(Q_{t-\tau_x} (P_sf)\ind_{\{\tau_x<t\}})+\bE(Q_{t+s-\tau_x}(f)\ind_{\{\tau_x\in[t,t+s]\}})\\
   =\, & \bP(\tau_x>t)\bP(\tau_x>s)  +\bE(Q_{t-\tau_x}(P_sf)\ind_{\{\tau_x<t\}})\\
&  + \textup{e}^{-t/\xi_x}\bE(Q_{s-\tau_x}(f)\ind_{\{\tau_x\in[0,s]\}})\\
 =\,  & \bP(\tau_x>t)\big(\bP(\tau_x>s)  +\bE(Q_{s-\tau_x}(f)\ind_{\{\tau_x\in[0,s]\}})\big)\\
& + \bE(Q_{t-\tau_x} (P_sf)\ind_{\{\tau_x<t\}})\\
 =\, & P_tP_sf(x), 
\end{align*}
which establishes the semigroup property of $P$. Finally, the distributional identity in \textit{(d)} is inherited from the $X^{(N)}$, because $\textup{tr}_{[N]}(\cdot)$ is continuous w.r.t.\ convergence in distribution.
\end{proof}
We interpret the environments $\xi$ as elements of the space $\mathcal{S}^{\alpha}_+\subset \ell^1$ given by
\[
\mathcal{S}^{\alpha}=\{\xi: \sup_{x\in\bN} |\xi_x| x^{1/\alpha}<\infty\},\; \mathcal{S}^{\alpha}_+=\mathcal{S}^\alpha \cap \{\xi\in\ell^1:\xi\geq 0 \}.
\]
Furthermore, let $\mathcal{M}(\bN)$ denote the space of probability measures on $\bN$ equipped with $\ell^1$-metric (which in our discrete setting is equivalent to the total variation metric) and let $ (\mathcal{B}(\bN), \|\cdot\|)$ be the space of bounded linear operators on $\bN$ with $\|\cdot\|$ denoting the canonical operator norm $\|L\|= \sup_{f:|f|_{\ell^\infty}\leq 1}|Lf|_{\ell^\infty}$. The next lemma addresses continuity of $P$ and $Q$ in the environment variable $\xi$.
\begin{lemma}\label{lemma:P-continuous}
Let $\alpha$ be fixed. For every $t>0$, the maps $Q_t: (\mathcal{S}^\alpha_+, \ell^1)\to (\mathcal{M}(\bN), \ell^1)$ and $P_t: (\mathcal{S}^\alpha_+, \ell^1)\to (\mathcal{B}(\bN), \|\cdot\|)$ are continuous.	
\end{lemma}
\begin{proof}
We begin by demonstrating the continuity of $Q_t$, for $\xi \in \mathcal{S}^{\alpha}_+$ set $q(\xi)=Q_t(\xi)$ for brevity. First note that the convergence of the characteristic functions established in Part (ii) of the proof of \ref{prop:Qexistence} implies that if $\xi^{(n)}\to \xi$ in $\ell^1$ with $\xi,\xi^{(n)}\in \mathcal{S}^{\alpha}_+$, then $q(\xi^{(n)})\to q(\xi)$ coordinate-wise, since the expression $\psi(\lambda)$ for the limiting characteristic function is continuous in $\xi$ w.r.t.\ the $\ell^1$-distance. It remains to extend the continuity to the stronger topology of $\ell^1$ convergence on $\mathcal{M}(\bN).$ 

To this end, let $r_M(\xi)$ denote the approximating measures $R^{M,\xi}$ used in the proof of Proposition~\ref{prop:Qexistence}. Fix $N$ and observe that the constants $C(\xi)$ in \eqref{eq:ConeandCtwo} can be chosen as a continuous function of the coordinate $\xi_{N+1}$ and the sums $\xi(\bN\setminus[N]),\xi([N])$. Note further that $N$ in \eqref{eq:ConeandCtwo} is chosen independently of $\xi$; consequently the bound there holds uniformly (with $C(\xi)$ replaced by some $C'$) in some small neighbourhood of $\xi$ w.r.t.\ the $\ell_1$-metric, i.e.\ for any $\epsilon>0, N\in\bN$, we can find some $\delta>0$ with
\[
r_M(\xi'; \bN\setminus[N])\leq \epsilon+ C'(\delta,\xi)tN^{1-1/\alpha},\quad M>\frac{N}{\epsilon}, \xi'\in \{\eta:|\xi-\eta|_{\ell_1}<\delta\}\equiv U(\delta,\xi).
\]
In the limit $M\to\infty$, the dependence of $M$ on $\epsilon$ vanishes (mind that in the calculation leading to \eqref{eq:ConeandCtwo}, $\epsilon$ was needed to deal with the initial condition in the finite system), which yields the uniform bound
\begin{equation}\label{eq:uniformQtail}
q(\xi';\bN\setminus[N])\leq C'(\delta,\xi)tN^{1-1/\alpha}, \xi'\in U(\delta,\xi).
\end{equation}
We conclude that
\begin{equation}\label{eq:q_unif}
|q(\xi)-q(\xi')|_{\ell^1}\leq \sum_{x=1}^N |q(\xi,\{x\})-q(\xi',\{x\})| + C'(\xi,\delta)tN^{1-1/\alpha}, \quad \xi'\in U(\delta,\xi),
\end{equation}
from which continuity of $q(\cdot)$ follows.

To address continuity for $P_t$ in $\xi$, we write $P_t=P_t(\xi), P'_t=P_t(\xi')$ and bound
\begin{equation}
	\begin{aligned}
	|P_tf-P'_tf|_\infty &\leq \sup_x |f(x)\bP(\tau_x>t)-f(x)\bP(\tau'_x>t)| \\
	& \phantom{=} + \sup_x\left|\int_0^t Q'_{t-s}f\frac{1}{\xi'_x}\textup{e}^{-s/\xi'_x}-Q_{t-s}f\frac{1}{\xi_x}\textup{e}^{-s/\xi_x} \di s \right|.  
	\end{aligned}
\end{equation}
Since $|f(x)|\leq 1$, the first term is at most $$\sup_{x}|\textup{e}^{-t/\xi_x}-\textup{e}^{-t/\xi'_x}|$$ which vanishes as $\xi'\to\xi$ in $\ell^1$. Similarly, the second term is bounded by
\begin{equation}
	\begin{aligned}
		 \sup_x \int_{0}^t |(Q'_{t-s}-Q_{t-s})| \frac{1}{\xi'_x}\textup{e}^{-s/\xi'_x}\di s + \sup_x |F_{\tau'(x)}(t)-F_{\tau(x)}(t)|, 
	\end{aligned}
\end{equation}
where $F_{\tau_x}$ denotes the distribution function of $\tau_x$. It is straightforward to see, using that $C(\xi',\delta)$ does not depend on $t$ in \eqref{eq:q_unif}, that both terms vanish for any $t>0$ as $\xi'\to\xi$ in $\ell^1$.
\end{proof}

The process $X$ enters $\bN$ according to the entrance law $(Q_t)_{t\geq0}$ at some site $x$ and the random walk dynamics correspond to waiting an $\operatorname{Exp}(\xi_x^{-1})$-distributed time at site $x$, then jumping to infinity and instantaneously re-entering $\bN$ again via $(Q_t)_{t>0}$. We next determine the invariant distribution of $X$. 
\begin{lemma}\label{lem:stationary}
The invariant distribution of $X$ is given by $\mu(x)=\frac{1}{\xi(\bN)} \xi_x, x\in\bN$.
\end{lemma}
\begin{proof}
$X^{(N)}$ is clearly ergodic, hence has a unique invariant distribution $\mu_N$. It is straightforward to deduce from the occupation time formula
$$
\lim_{t\to\infty}\frac{\int_0^t \ind_{\{X^{(N)}_s=x\}}\di s}{t}=\mu_N(x),\; \bP\text{-almost surely for } x\in\bN,
$$
 that $\mu_N(x)=\xi_x/\xi([N])$, $x\in[N]$, is the invariant distribution of $X^{N}$. According to Theorem~\ref{thm:OSTRW} (c) the semigroups $P^{(N)}$ converge to $P$ and we clearly have $\lim_{N\to\infty}\mu_N=\mu$, hence \cite[Theorem 9.10]{EthieKurtz86} or a direct calculation yield that $\mu$ is invariant for $X$.
\end{proof}

\begin{lemma}
Let $f,g:\bN \to \bR$ have finite support. Then 
\[ \lim_{N\to\infty} \langle P_t^{(N)} f, g\rangle_{\mu_N} = \langle P_t f, g\rangle_{\mu},\]
and 
\[\langle P_t f, g\rangle_{\mu} = \langle f, P_t g\rangle_{\mu}.\]
\end{lemma}
\begin{proof}
The convergence result is a straightforward consequence of the convergence of $\mu_N$ to $\mu$ and of $P^{(N)}$ to $P$. Self-adjointess of $P$ w.r.t.\ $\mu$ is inherited from $P^{(N)}$, since $X^{(N)}$ is reversible.
\end{proof}
\begin{lemma}
Let $p_t^{(N)}(x,y) = \bP_x(X_t^{(N)}=y)$. Then, for $x\neq y$, 
$$ \lim_{N\to\infty}\frac{\di}{\di t}p_t^{(N)}(x,y)|_{t=0}= 0, $$
whereas $$\lim_{N\to\infty}p_t^{(N)}(x,y) > 0, \quad \text{for all }t>0.$$
\end{lemma}
\begin{proof}
The first statement is simply a reformulation of the fact that the transition rates $r^{(N)}(x,y)$ vanish as $N\to\infty$. To prove the second statement we rewrite$$
p_t(x,y)=\delta_xP_t\ind_{\{y\}} = \bE(Q_{t-\tau_x}(y)\ind_{\{\tau_x<t\}}).
$$
The claim follows once we show $\supp Q_s=\bN$ for any $s>0$. Note that for every $x\in\bN$ and $t>0$ there exists some $K>x$ such that $Q_{t}([K])>0$, hence there is a positive probability that $X$ has visited $[K]$ before time $t$. Upon entering $[K]$ at time $T_0<t$, the position $X_{T_0}\in[K]$ is uniformly distributed by Theorem~\ref{thm:OSTRW} (d), and thus there is a positive probability that $X_s=x$ for all $s\in[T_0,t]$, hence $Q_t(x)>0$.
\end{proof}
\begin{lemma}\label{lem:coupling2RW}For any $x,y\in \bN$ we have
\begin{align}
\norm{\bP_x(X_t \in \cdot) - \bP_y(X_t \in \cdot)}_{\textup{TV}} \leq e^{-\xi_{1}^{-1}t},
\end{align}
and in particular
\begin{align}
\norm{\bP_x(X_t \in \cdot) - \mu}_{\textup{TV}} \leq e^{-\xi_1^{-1}t}.
\end{align}
\end{lemma}
\begin{proof}
We construct a Markovian coupling of two versions of the limit random walk $X_t$ and $Y_t$. If $X_t<Y_t$, then the exit rate $\xi_{Y_t}^{-1}$ from the current state of $Y_t$ is larger than the corresponding rate for $X_t$. Hence, we can couple the two exponentials governing the exit so that when $X_t$ jumps, $Y_t$ jumps as well. In the case of a joint exit, we can couple the random walks to move identically, as both movements are prescribed by $(Q_t)_{t>0}$. The additional jumps of $Y_t$ may cause $Y_t$ to jump to the left of $X_t$, in which case the roles of $X_t$ and $Y_t$ reverse so that it is always the random walk with the lower index, which determines the coupling speed. Since $\xi_{1}$ is the largest weight, the exit time for a random walk started in the slowest site dominates this coupling argument, which proves that the coupling time in this coupling is at most an exponential with parameter $\xi_1^{-1}$.\end{proof}
Note that the coupling used in the proof of Lemma~\ref{lem:coupling2RW} is not optimal, as it ignores the possibility that $Y_t$ jumps on top of $X_t$ for an earlier coupling. We proceed by stating two direct implications of the lemma.
\begin{corollary}\label{Cor:spec}
We have $\bfP$-almost surely that 
\begin{enumerate}[(a)]
\item $\norm{Q_t-\mu}_{\textup{TV}}\leq e^{-\xi_1^{-1}t};$
\item $(P_t)$ has a spectral gap of at least $\xi_1^{-1}$.
\end{enumerate}
\end{corollary}

To construct a system of infinitely many coalescing versions of $X$ in Section 4, it is technically convenient to extend the trajectory space.
\begin{lemma}\label{lem:Feller}
Let $\bar{\bN}$ denote the one-point compactification of $\bN$ with respect to the discrete topology. Then $X$ can be extended to a Feller process on $\bar{\bN}$.
\end{lemma}
\begin{proof}
We do not distinguish notationally between $X$ and its extension. We use the probabilistic characterisation of the Feller property given in \cite[Lemma 17.3]{KallenbergFM}. This amounts to establishing the two conditions
\begin{equation}\label{eq:F1}
\bP_x(X_t \in \cdot)\to \bP_y(X_t \in \cdot) \text{ for any }t\geq 0 \text{ as }x\to y \text{ in }\bar{\bN},
\end{equation}
and
\begin{equation}\label{eq:F2}
X_t\overset{\bP}{\to}X_0 \text{ as }t\to 0.
\end{equation}
If $y\in \bN$, then \eqref{eq:F1} holds trivially. If $y=\infty$, then as $x\to y$ the time of the first jump converges to $0$ in probability and thus $P_t\delta_x$ converges to $Q_t$ as required for \eqref{eq:F1}. Similarly, item \eqref{eq:F2} is trivial if the process is initiated in $\bN$. If the process is initiated via the entrance law, i.e.\ it starts at $\infty$, then $X_t$ must leave any finite subset of $\bN$ as $t\to 0$. But this is equivalent to $X_t\to\infty$ in $\bar{\bN}$, which yields \eqref{eq:F2}.
\end{proof}


\subsection*{Multiple walks.}
We now consider systems of finitely many independent limit random walks. Let $\bP_{\nu_1\otimes\nu_2}$ denote the distribution of two independent copies $(X_t)_{t\geq 0}$ and $(Y_t)_{t\geq 0}$ of the limiting walk $X$, started in distributions $\nu_1$ and $\nu_2$, respectively. If the starting positions are fixed $\nu_1=\delta_x, \nu_2=\delta_y$, then we write $\bP_{(x,y)}$ instead. We say that $X$ and $Y$ \emph{meet} at time $t$, if $X_t=Y_t$.
\begin{corollary}
The probability that two independent random walks with fixed starting points meet at time $t$ satisfies
\begin{align}
    \bP_{(x,y)}(X_t = Y_t) \geq \sup_{\alpha\in(0,1)}(1-\alpha)^2\sum_{z\leq h_\alpha(t)}\mu(z)^2,
\end{align}
where $h_\alpha(t) = \sup\{x\in\bN : \mu(x) \geq \alpha^{-1}e^{\xi_1^{-1}t} \}$.
\end{corollary}
\begin{proof}
Fix $\alpha\in(0,1)$, then we have
\begin{align}\label{eq:mult1}
\bP_{(x,y)}(X_t = Y_t) &= \sum_z p_t(x,z)p_t(y,z)   \\
&= \sum_{z\leq h_\alpha(t)} \big(\mu(z) + (p_t(x,z)-\mu(z)) \big) \big(\mu(z) + (p_t(y,z)-\mu(z)) \big)      \\
&\geq \sum_{z\leq h_\alpha(t)} (\mu(z) - \norm{p_t(x,\cdot)-\mu}_{\textup{TV}} ) (\mu(z) - \norm{p_t(y,\cdot)-\mu}_{\textup{TV}} )\\
&\geq \sum_{z\leq h_\alpha(t)} (\mu(z) - e^{-\xi_1^{-1}t} ) (\mu(z) - e^{-\xi_1^{-1}t} ),
\end{align}
where we used Corollary \ref{Cor:spec} in the last inequality and the fact that $\mu(z)>\norm{p_t(x,\cdot)-\mu}_{\textup{TV}}> \sup_{w\in\mathbb{N}} |p_t(x,w)-\mu(w)|$, for $z\in h_\alpha(t)$  in the second inequality.
Note, that by definition of $h_\alpha(t)$, the last expression is at least $(1-\alpha)^2\sum_{z\leq h_\alpha(t)}\mu(z)^2$. The claim now follows by taking the supremum in $\alpha$.
\end{proof}
It follows that for $t$ of order $\xi_1^{-1}$ there is a realistic chance that two independent random walks have met. This time-scale is also the best possible if we make no assumptions on the starting points, since on time scales much smaller than $\xi_1^{-1}$, two random walks starting in $1$ and $2$ have not jumped with high probability, hence they cannot have met. We next show that if we start both walks in the stationary distribution, they typically do not meet in high index sites. 

\begin{lemma}\label{lem:nomeet}
There is some constant $C<\infty$ (independent of $\xi$) such that for any $A\subset\bN$ and $T>0$, $\bfP$-almost surely
$$\bP_{\mu\otimes\mu}(X \text{ and } Y \text{ meet inside } A \text{ before time }T)\leq C\Big(T \,\frac{\mu(A)}{\xi(\bN)} + \sum_{z\in A}\mu(z)^2\Big).$$
\end{lemma}
\begin{proof}
Let  $$\tau_z=\inf\{t: X_t=Y_t=z\}, z\in \bN.$$ 	
Since $X$ and $Y$ individually leave $z$ at rate $\xi_z^{-1}$, the joint leaving rate is $2\xi_z^{-1}$. Let $\sigma_z$ be the waiting time for the next jump of either walker after $\tau_z$. Note that $\sigma_z$ is independent of $\tau_z$ by the Markov property. Then
\begin{align}
&\bP_{\mu\otimes\mu}\left(\tau_z\leq T \right) 
= \bP_{\mu\otimes\mu}\left(\tau_z\leq T, \sigma_z\geq \epsilon \xi_z \right)      
+ \bP_{\mu\otimes\mu}\left(\tau_z\leq T, \sigma_z< \epsilon \xi_z \right)
\\
&\leq\bP_{\mu\otimes\mu}\left(\int_0^{T+\epsilon\xi_z} \ind_{\{X_t=z\}}\ind_{\{Y_t=z\}}\di t \geq \epsilon \xi_z \right) 
+ \bP_{\mu\otimes\mu}\left(\tau_z\leq T\right)\bP\left(\sigma_z< \epsilon \xi_z \right).
\end{align}
Note that the increase in the time interval from $[0,T]$ to $[0,T+\epsilon\xi_z]$ accommodates for the event that $\tau_z\in(T-\epsilon\xi_z,T]$. Since $\bP\left(\sigma_z< \epsilon \xi_z \right)=1-e^{-2\epsilon}$, we obtain by regrouping terms that
$$
\bP_{\mu\otimes\mu}\left(\tau_z\leq T \right)\leq e^{2\epsilon} \bP_{\mu\otimes\mu}\left(\int_0^{T+\epsilon\xi_z} \ind_{\{X_t=z\}}\ind_{\{Y_t=z\}}\di t \geq \epsilon \xi_z\right).
$$
By Markov's inequality,
\begin{align}
\bP_{\mu\otimes\mu}\left(\int_0^{T+\epsilon\xi_z} \ind_{\{X_t=z\}}\ind_{\{Y_t=z\}}\di t \geq \epsilon \xi_z \right)       
&\leq \frac{\bE_{\mu\otimes\mu} \int_0^{T+\epsilon\xi_z} \ind_{\{X_t=z\}}\ind_{\{Y_t=z\}}\di t }{\epsilon \xi_z}   \\
&= \frac{(T+\epsilon\xi_z)\mu(z)^2}{\epsilon \xi_z} = \Big(\xi_z+\frac{T}{\epsilon}\Big) \frac{\xi_z}{\xi(\bN)^2}.
\end{align}
Therefore, 
\begin{align}\label{eq:meeting-estimate-1}
\bP_{\mu\otimes\mu}\left(\tau_z\leq T \right)
\leq e^{2\epsilon} \Big(\xi_z+\frac{T}{\epsilon}\Big) \frac{\xi_z}{\xi(\bN)^2}
\end{align}
Choosing $\epsilon$ of order one and applying a union bound over all $z\in A$ gives the desired result.
\end{proof}

For arbitrary starting positions, we obtain the following corollary by using the simple fact that $\bP_{(x,y)}(E)\leq \frac{1}{\mu(x)\mu(y)} {\bP_{\mu\otimes\mu}}(E)$ for any event $E$.

\begin{corollary}\label{cor:meetingbound}
Let $x<y$ and $A\subset \bN$, then
$$\bP_{(x,y)}(X \text{ and } Y \text{ meet in } A \text{ before time }T)\leq \frac{C}{\xi_x\xi_y}\big(T \,\xi(A)  +{\sum_{z\in A}\xi_z^2}\big).$$
\end{corollary}

%
%
%

\section{From ranked environments to the limit coalescing walk}
The limiting walk $X$ of Section~3 was constructed from the environment $\xi$. We now isolate the abstract convergence statement that connects the finite models to this limit. Given positive weights $w_1^{(N)}\geq\dots\geq w_N^{(N)}$ and a scale $a_N$, define the rescaled ranked environment
\[
\widehat{w}^{(N)}=\big(a_N^{-1}w_1^{(N)},\dots,a_N^{-1}w_N^{(N)},0,0,\dots\big)\in \ell_1.
\]
By the $\ell_1$-criterion established at the beginning of the previous section, convergence $\widehat{w}^{(N)}\overset{d}{\to}\xi$ in $\ell_1$ is equivalent to pointwise convergence of the ranked weights together with convergence of the total mass. In the heavy-tailed setting introduced in Section~1, this is precisely the extreme-value input \eqref{eq:limweights}. The results of the present section use only this ranked-environment convergence and no further properties of the microscopic law.

\subsection*{Single-lineage convergence}
\begin{proposition}[Abstract ranked-environment convergence]\label{prop:convYtoX}
	Let $w^{(N)}_1\geq w^{(N)}_2\geq\dots\geq w^{(N)}_N$ be positive weights and let $\widehat{w}^{(N)}$ be as above. Assume that
	\[
	\widehat{w}^{(N)}\overset{d}{\longrightarrow}\xi \qquad \text{in }\ell_1,
	\]
	where $\xi$ is the limiting environment used to construct $X$ in Section~3. Let $Y^{(N)}$ denote the random walk in the environment $(w^{(N)}_i)_{i\in[N]}$ started in position $y_0\leq M$ with uniform jump distribution. Then
	\[
	Y^{(N)}_{a_N(\cdot)}\overset{d}{\longrightarrow} X_{(\cdot)},\quad\mbox{as }N\to\infty.
	\]
\end{proposition}
\begin{proof}
	It suffices to show that $\big(Y^{(N)}_{a_N t}\big)_{t\in[0,T]}$ can be perfectly coupled to $X^{(N)}$ (which coincides with the trace of $X$ on $[N]$ by Theorem~\ref{thm:OSTRW}(d)) up to some error probability that vanishes as $N\to\infty$. Note that given the initial positions $Y_0^{(N)},\;X_0$ the processes can be reconstructed from the environment of exponential clocks. Hence, it suffices to show that the Poisson field of clock processes corresponding to $w^{(N)}$ at time scale $a_N$ converges weakly to the Poisson field of clock processes parametrised by $\xi_1,\dots,\xi_N$. The assumed convergence $\widehat{w}^{(N)}\overset{d}{\to}\xi$ in $\ell_1$ implies pointwise convergence of the clock intensities together with convergence of their total mass, and therefore convergence of the corresponding characteristic functionals. The claim now follows from the fact that a sequence of Poisson measures converges weakly to a limiting Poisson measure if and only if the corresponding characteristic functionals converge, see e.g.\cite[Proposition 11.1.VIII]{DVJII}.
\end{proof}
In particular, once the ranked environments converge, the motion of a single lineage is universal on the time scale $a_N$. For the complete-graph models of interest here, the domain-of-attraction statement from the introduction is obtained by combining Proposition~\ref{prop:convYtoX} with \eqref{eq:limweights}.

\subsection*{Finite coalescing systems}
To pass from the motion of one lineage to the voter model via duality, we next extend the convergence to finitely many coalescing lineages. The only new issue is to control coalescences that take place far out among the fast sites. Corollary~\ref{cor:meetingbound} provides precisely the required estimate.
Let us now formally define the coalescence operator $\operatorname{coal}(\cdot)$ for an arbitrary collection $(X^{1},X^{2},\cdots)$ of copies of $X$ with distinct initial conditions. We assume $X^{1}_0=1, X^{2}_0=2,\dots$. This induces a natural order on the paths: $X^{i}<X^{j}$ if and only if $X^{i}_0<X^{j}_0$. Note that this convention imposes no additional restriction since an arbitrary collection of paths can always be ordered according to starting positions and viewed as a projection of the full configuration. For any $i$ denote by $$\sigma_{i}=\inf\{t: X_t^{i}=X_t^{j} \text{ for any }j<i\}$$ the \emph{coalescence time} of the $i$-th walker. Note that each $\sigma_i$ is measurable with respect to the $\sigma$-field $\sigma(X^{j}, j\leq i)$ which is generated by only finitely many walkers and that $\sigma_i>0$ almost surely for each $i$. In particular the coalescence time of any walker is well-defined and there are only countably many coalescence times. Define now
\begin{equation}\label{eq:coaldefi}
\operatorname{coal}\big((X^{1},X^{2},\dots)\big)_t=\Big(\sum_{j: \sigma_j> t}\delta_{X^{j}_t}(x)\Big)\in \{0,1\}^{\bN}, \quad t\geq 0.
\end{equation}
The existence of this infinite system of coalescing walks is deduced from standard results in the next proposition.
\begin{proposition}\label{prop:projlim}
For each finite $M\subset\bN$, let $\mathbb{Q}_M$ denote the distribution of
\[
\operatorname{coal}\big((X^{j}, j\in M)\big).
\]
Then
$$
\mathbb{Q}=\mathcal{L}\big(\operatorname{coal}\big((X^{1},X^{2},\dots)\big)\big)
$$
is the projective limit of $\{\mathbb{Q}_M, M\subset\bN \text{ finite}\}.$
\end{proposition}
\begin{proof}
It is clear that if the projective limit exists, it must coincide with $\mathbb{Q}$ and this would also establish the well-definedness of \eqref{eq:coaldefi}. To this end, we invoke a version of the Kolmogorov extension theorem, namely \cite[Theorem 8.23]{KallenbergFM}. Recall that by Lemma~\ref{lem:Feller} we may view the particle trajectories as paths in $\bar{\bN}$, which is Borel isomorphic to $\{1/n, n\in\bN\}\cup\{\infty\}$ and thus a Borel space. Hence, we only need to show the consistency of the family $\{\mathbb{Q}_M, M\subset\bN \text{ finite}\}$. The latter is an immediate consequence of the definition of the coalescence times $\sigma_j, j\in M$, if for each particle, one keeps track of the indices of the particles that have already merged with it\footnote{ Formally, this requires to extend the probability spaces involved. However, we only use the additional structure in this proof and we found that the corresponding notation only obscures the argument}: Arguing inductively, it suffices to show that $\mathbb{Q}_{M'}$ coincides with the projection of $\mathbb{Q}_M$ onto the particles with indices in $M'$ for $|M'|=|M|-1$. Denote the index of the additional (\emph{spare}) particle in $M$ by $s$. If no other particle coalesces with $s$ before $\sigma_s$, then there is nothing to show. Assume this is not the case, then there exist coalescence times $\sigma_{j_1}<\sigma_{j_{2}}<\dots<\sigma_{j_k}<\sigma_s$ at which particles with indices in $M'$ have merged with the spare particle. Since these \emph{virtual} coalescences do not occur in the smaller system, we need to define the trajectories of the non-spare particles involved after the virtual coalescence events. By coupling the $M'$ to the $M$ system, we may declare the trajectory of the particle $j_1$ involved in the first virtual coalescence to be that of the spare particle until the next coalescence event involving $j_1$. By construction, the latter will occur at time $\sigma_{j_{2}}$. Depending on whether $j_1<j_2$ or $j_1>j_2$, one of the particles survives and we set it once again to follow the trajectory of the spare particle. We continue this procedure and thus inductively define a path for the $M'$-coalescing random walk. It follows from the Markov property of the underlying walks that this path has indeed the distribution $\mathbb{Q}_{M'}$.
\end{proof}
Note that the law $\mathbb{Q}$ just constructed formally corresponds to a coalescing system with trajectories in $\{0,1\}^{\bN}\times \bN$, since by convention we do not coalesce particles that meet at $\infty$. That this is justified follows from Corollary~\ref{cor:meetingbound}: the probability that particles meet `near' the boundary point $\infty$ vanishes in the finite system, which is also crucially used in the proof of Proposition~\ref{prop:convYtoXfinite} below. However, we may also consider the process $\gamma$ derived from propagating an initial condition $\gamma_0$ according to (projections of) $\mathbb{Q}$ as a process on $\{0,1\}^{\bN}$, just as in Sections 1 and 2.

\begin{proposition}[Finite-lineage convergence]\label{prop:convYtoXfinite}
Let $\gamma_0^{(N)}\overset{d}{\to}\gamma_0\in\{0,1\}^{\bN}$ as $N\to\infty$, where $\supp \gamma_0^{(N)}\subset[K]$ for some $K$, uniformly in $N\in\bN$. Then the corresponding coalescing random walk systems 
\[
\big(\gamma_{a_N t}^{(N)}\big)_{t\geq 0}
\]
admit a distributional limit $\gamma=(\gamma_{t})_{t\geq 0}$ as $N\to\infty$, in which individual particle dynamics are given by i.i.d.\ copies of $X$ and particles coalesce whenever they share the same location.
\end{proposition}
\begin{proof}
	We define the following approximations to the coalescence operator:
	\begin{itemize}
		\item \(\operatorname{coal}_M(\cdot)\) denotes the coalescence operator allowing coalescences only at sites in $[M]$;
		\item \(\operatorname{coal}_{M,\epsilon}(\cdot)\) denotes the coalescence operator allowing coalescences only at sites in $[M]$ and at times in $\epsilon\bN$.
	\end{itemize}
	Set
	\[
	\mathbf{Y}^{(N)}=(Y^{(1,N)},\dots,Y^{(k,N)}),\qquad \mathbf{X}=(X^{(1)},\dots,X^{(k)}).
	\]
	By Proposition~\ref{prop:convYtoX}, we have that for any $M\in\mathbb{N}$, $\epsilon>0$, $\vec{t}=(t_1,\dots,t_d)\in[0,T]^d$ and $\vec{x}\in \bigtimes_{j=1}^d \{0,1\}^{\bN}$
	\begin{equation}\label{eq:limitprobab1}
		\begin{aligned}
		\bP\left(\operatorname{coal}_{M,\epsilon a_N}(\mathbf{Y}^{(N)})_{a_N\vec{t}}=\vec{x}\right)&\\
		 \overset{N\to\infty}{\longrightarrow}\bP\Big(&\operatorname{coal}_{M,\epsilon}(\mathbf{X})_{\vec{t}}=\vec{x}\Big).
		\end{aligned}
	\end{equation}
	We next remove the temporal discretisation and the spatial truncation on both sides of \eqref{eq:limitprobab1}.

	First consider the finite system. Under $\operatorname{coal}_{M}$ there are at most $k-1$ coalescence events. If one of them is not detected by $\operatorname{coal}_{M,\epsilon a_N}$, then after the corresponding coalescence event at some site in $[M]$ the merged particle must jump within the following $\epsilon a_N$ time units. Hence
	\begin{equation}\label{eq:finite-time-disc}
		\bP\Big( \operatorname{coal}_{M,\epsilon a_N}(\mathbf{Y}^{(N)})_{a_N\vec{t}}\neq \operatorname{coal}_{M}(\mathbf{Y}^{(N)})_{a_N\vec{t}} \Big)
		\leq (k-1)\Big(1-\textup{e}^{-\epsilon a_N/w_M^{(N)}}\Big).
	\end{equation}
	By \eqref{eq:limweights}, we have $a_N/w_M^{(N)}\to \xi_M^{-1}$ for each fixed $M$, so that for fixed $M$
	\begin{equation}\label{eq:finite-time-disc-limit}
		\limsup_{N\to\infty}\bP\Big( \operatorname{coal}_{M,\epsilon a_N}(\mathbf{Y}^{(N)})_{a_N\vec{t}}\neq \operatorname{coal}_{M}(\mathbf{Y}^{(N)})_{a_N\vec{t}} \Big)
		\leq (k-1)\Big(1-\textup{e}^{-\epsilon/\xi_M}\Big),
	\end{equation}
	which vanishes as $\epsilon\downarrow 0$.

	Next, the error caused by suppressing coalescences outside $[M]$ is bounded by
	\begin{align}
		& \bP\big( \operatorname{coal}_{M} (\mathbf{Y}^{(N)})_{a_N\vec{t}}\neq \operatorname{coal}(\mathbf{Y}^{(N)})_{a_N\vec{t}} \big)\nonumber\\
		&\leq \binom{k}{2} \sup_{x,y\in [K]} \bP_{(x,y)}\Big( Y^{1,N}\text{ and } Y^{2,N} \text{ meet before time }Ta_N \text{ in }[N]\setminus[M]\Big). \label{eq:finite-space-disc}
	\end{align}
	By the same argument that led to Corollary~\ref{cor:meetingbound},
	\[
	\begin{aligned}
	&\bP_{(x,y)}\Big( Y^{1,N}\text{ and } Y^{2,N} \text{ meet before time }Ta_N \text{ in }[N]\setminus[M]\Big)\\
	&\qquad =O\left(\frac{T a_N \sum_{z=M+1}^N w^{(N)}_z}{w^{(N)}_x w^{(N)}_y} \right),
	\end{aligned}
	\]
	uniformly in $x,y\in[K]$. Moreover, for each fixed $x\in[K]$ we have
	$a_N^{-1}w_x^{(N)}\to\xi_x$, and the tail masses satisfy
	$a_N^{-1}\sum_{z=M+1}^N w^{(N)}_z\to \xi(\bN\setminus[M])$.
	It follows that for every fixed $M$
	\begin{equation}\label{eq:finite-space-disc-limit}
		\limsup_{N\to\infty}\bP\big( \operatorname{coal}_{M} (\mathbf{Y}^{(N)})_{a_N\vec{t}}\neq \operatorname{coal}(\mathbf{Y}^{(N)})_{a_N\vec{t}} \big)
		\leq C(T,K,\xi)\,\xi(\bN\setminus[M]),
	\end{equation}
	and the right hand side vanishes as $M\to\infty$.

	The same two estimates hold for the limiting system. Indeed, since after a coalescence event in $[M]$ the merged limit particle leaves its current site at rate at most $\xi_M^{-1}$, we have
	\begin{equation}\label{eq:limit-time-disc}
		\bP\Big( \operatorname{coal}_{M,\epsilon}(\mathbf{X})_{\vec{t}}\neq \operatorname{coal}_{M}(\mathbf{X})_{\vec{t}} \Big)
		\leq (k-1)\Big(1-\textup{e}^{-\epsilon/\xi_M}\Big),
	\end{equation}
	which tends to $0$ as $\epsilon\downarrow 0$. Likewise,
	\begin{align}
		& \bP\big( \operatorname{coal}_{M} (\mathbf{X})_{\vec{t}}\neq \operatorname{coal}(\mathbf{X})_{\vec{t}} \big)\nonumber\\
		&\leq \binom{k}{2}\sup_{x,y\in[K]}\bP_{(x,y)}\Big( X^{1}\text{ and } X^{2} \text{ meet before time }T \text{ in }\bN\setminus[M]\Big)\nonumber\\
		&\leq C(T,K,\xi)\Big(\xi(\bN\setminus[M])+\sum_{z>M}\xi_z^2\Big), \label{eq:limit-space-disc}
	\end{align}
	by Corollary~\ref{cor:meetingbound}, and the right hand side vanishes as $M\to\infty$.

	Combining \eqref{eq:limitprobab1}, \eqref{eq:finite-time-disc-limit}, \eqref{eq:finite-space-disc-limit}, \eqref{eq:limit-time-disc} and \eqref{eq:limit-space-disc}, we may first let $N\to\infty$, then $\epsilon\downarrow 0$ and finally $M\to\infty$. This shows that the finite-dimensional distributions of $\big(\gamma_{a_N t}^{(N)}\big)_{t\geq 0}$ converge to those of the coalescing system generated by i.i.d.\ copies of $X$.
\end{proof}

\section{Coalescence time and coming down from infinity}
In this section, we establish that $\gamma$ `comes down from infinity'.
\begin{theorem}\label{thm:comingdown}
Let $\gamma_0\equiv 1$. For almost every environment $\xi$, we have that
\[
\sum_{x\in\bN}\gamma_t(x)<\infty \quad \text{ for any }t>0.
\]
Furthermore, $\gamma$ can be extended to a Feller process on $\{0,1\}^{\bN}\times 2^{\bN}$ that enters $\{0,1\}^{\bN}$ continuously.
\end{theorem}
Theorem~\ref{thm:comingdown} is proved at the end of this section. The main ingredient for the proof is the following auxiliary statement in which we control the expected coalescence time $\bE \tau_{M+1\to M}$. Here, $\tau_{k\rightarrow l}$ denotes the first time at which the process initiated with $k$ particles enters a state with $l$ particles for $l<k$. In fact, it suffices to consider the process $\gamma$ restricted to sites in $[N].$
\begin{proposition}\label{thm:coalescence-time}
	There exists a constant $C(\xi)$ so that, for any $M$,
	\begin{align}
		\bE \big[\tau^\super{N}_{M+1\to M}\big | \gamma_0^{(N)}\big] \leq C(\xi)M^{-\frac1\alpha},
	\end{align}
	uniformly in $N$ and $\gamma_0^{(N)}\in\{0,1\}^{[N]}$.
\end{proposition}

\subsection*{The proof of Proposition \ref{thm:coalescence-time}}
The proof is a direct consequence of Lemmas~\ref{lem:5.1}--\ref{lem:5.3} in this paragraph. Since a coalescence can only happen when a particle jumps we need to control the total number of jumps $J^\super{N}_t$ in a system of $M+1$ coalescing random walks on $[N]$, which we do by means of a lower bound on $J_t^\super{N}$.
\begin{lemma}\label{lem:5.1}
	Assume there is a lower bound $J'_t \leq J^\super{N}_t$ valid until $\tau^\super{N}_{M+1\to M}$, i.e.\ $J'_t \leq J^\super{N}_t + \infty\ind_{t> \tau\super{N}_{M+1\to M}}$ uniformly in all starting configurations $\gamma_0^{(N)}\in\{0,1\}^{[N]}$. Assume further that there are constants $c_1=c_1(\xi)$, $c_2=c_2(\xi)$ satisfying
	\begin{align}
		\int_0^t \bP\Big(J'_t < c_1 N M^{\frac1\alpha-1} t \Big) \di t \leq c_2 M^{-\frac1\alpha}.
	\end{align}
	Then
	\begin{align}
		\bE \big[\tau^\super{N}_{M+1\to M}\big | \gamma_0^{(N)}\big] \leq (c_1^{-1}+c_2)M^{-\frac1\alpha}.
	\end{align}
\end{lemma}
\begin{proof}
	Each jump in the system has probability $\frac{M}{N}$ to land on top of another particle, leading to a coalescence. Since these jumps are independent, the number of jumps needed is a geometric random variable $Z$ with success parameter $\frac{M}{N}$, and $\{\tau^\super{N}_{M+1\to M}>t\} = \{Z>J^\super{N}_t\}$. Note that $Z$ and $J^\super{N}_t$ are not independent. By assumption there is a coupling $\hP_\bfx$ with $J_t' \leq J_t^\super{N}$ on the event $\{\tau^\super{N}_{M+1\to M}>t\}$, where we set $\bfx=\gamma_0^{(N)}$ for conciseness. Then, writing $\tau = \tau^\super{N}_{M+1\to M}$, for any $a\in\bN$, 
	\begin{align}
		&\bP_\bfx \left( \tau>t \right) = \bP_\bfx\left( \tau>t, Z \geq a \right) + \bP_\bfx\left( \tau>t, Z<a\right)    \\
		&\leq \bP_\bfx\left( Z \geq a \right) + \bP_\bfx\left( \tau>t, J^\super{N}_t<a\right)    \\
		&\leq \bP_\bfx\left( Z \geq a \right) + \hP_\bfx\left( \tau>t, J'_t<a\right)    \\
		&\leq \bP_\bfx\left( Z \geq a \right) + \hP_\bfx\left(J'_t<a\right).
	\end{align}
	With $a=\lfloor c_1 NM^{\frac1\alpha -1}t \rfloor$,
	\begin{align}
		\bP_\bfx\left( Z \geq a \right) = \left(1-\frac{M}{N}\right)^{a-1}
		\leq e^{-c_1 M^{\frac1\alpha}t}
	\end{align}
	and 
	\begin{align}
		\int_0^\infty \bP_\bfx\left( Z \geq a \right) \;dt \leq c_1^{-1}M^{-\frac1\alpha}.
	\end{align}
	By assumption, 
	\begin{align}
		\int_0^t \hP_\bfx\left(J'_t<a\right) \;dt \leq c_2 M^{-\frac1\alpha}
	\end{align}
	and since $\bP_\bfx (\tau^\super{N}_{M+1\to M} =t )=0$
	\begin{align}
		&\bE_\bfx \tau^\super{N}_{M+1\to M} = \int_0^\infty \bP_\bfx \left(\tau^\super{N}_{M+1\to M} >t \right)\;dt \leq (c_1^{-1}+c_2)M^{-\frac1\alpha}. \qedhere
	\end{align}
\end{proof}

We now need to construct a good lower bound $J_t'$ on $J_t^\super{N}$. Let $X_t'$ be a random walk on $\{M+1,...,N\}$ started in $M+1$. The mean waiting times in site $x$ are given by $\xi_x$ (so identical to $X_t^\super{N}$), but it then jumps to $y$ with probability $\frac1N$ for $y>M+1$ and $\frac{M+1}{N}$ for $y=M+1$. Let $J_t'$ denote the number of jumps this random walk performs until time $t$.
\begin{lemma}
	The process $J'_t$ is a lower bound on $J^\super{N}_t$ until $\tau^\super{N}_{M+1\to M}$, that is, $J'_t \prec J^\super{N}_t + \infty\ind_{t> \tau^\super{N}_{M+1\to M}}$ uniformly in all starting configurations $\gamma_0^{(N)}$, where `$\prec$' indicates stochastic domination.
\end{lemma}
\begin{proof}
	Consider the system of coalescing random walks where we mark one particle, with $X_t''$ the position of the marked particle. Then, clearly, the number of jumps performed by $X''_t$ is a lower bound on the total number of jumps (we do not care what happens with the marked particle during a coalescence event, since the lower bound needs only to be valid until the first coalescence). The mechanism for marking is as follows: Initially, the particle with the largest position is marked, which is automatically at a position larger than $M$. Then, when the marked particle jumps to a site in $[M]$, we change the mark to the particle on the largest position, which is larger than $M$ (again, assuming coalescence hasn't happened yet). Note that by construction, until the first coalescence, $X''_t\geq M+1$.
	
	We now couple $X''_t$ with $X'_t$. If $X''_t = X'_t$, both particles move identically until $X''_t$ jump into $[M]$. On this event, the marked particle changes to the largest particle, so $X''_t \geq M+1$, with the position random depending on the position of the other particles. On the other hand, $X'_t$ jumps to $M+1$. If $X'_t = M+1$ and $X''_t>M+1$ the jump rate of $X_t''$ is $\xi_{X''_t}^{-1}$ and hence dominates the jump rate $\xi_{M+1}^{-1}$ of $X_t'$. Therefore we can couple the jump times so that if $X_t'$ jumps, $X''_t$ jumps as well and both jump to the same target (unless the target is in $[M]$, in which case $X_t'$ jumps to $M+1$ (and effectively does not move, but we do count this jump), and $X_t''$ jumps to some site $\geq M+1$.
	
	From this construction we obtain that whenever $X_t'$ jumps, so does $X_t''$, and hence $J_t'\leq J_t^\super{N}$ until the first coalescence time.
\end{proof}

\begin{lemma}\label{lem:5.3}
	$J'_t$ satisfies 
	\begin{align}
		\int_0^\infty \bP\Big(J'_t < \frac{N M^{\frac1\alpha-1} }{4(1+C(\xi))} t \Big)\;dt 
		\leq \frac4{C(\xi)}M^{-\frac1\alpha}.
	\end{align}
\end{lemma}
\begin{proof}
	Observe by Chebyshevs inequality
	\begin{equation}\label{eq:cheby}
	\bP(J'_t < a ) = \bP( \sum_{i=1}^a \sigma_i > t) \leq e^{-\beta t} \bE e^{ \beta\sum_i \sigma_i},
	\end{equation}
where $\sigma_i$ is the waiting time for the $i$-th jump. These are $\operatorname{Exp}(\xi_{X'_i}^{-1})$-distributed and independent given $X'_i$, $i=1,2,...$, where $X'_i$ is the position of $X'_t$ before the $i$-th jump. We have $X'_1=M+1$ and the remaining positions are i.i.d.\ with $\bP(X'_i=x) = \frac1N$, $x>M+1$ and $\bP(X'_i=M+1) = \frac{M+1}N$.
	Therefore, for $\beta\leq \frac12\xi_{M+1}^{-1}$ (and hence $\beta \xi_x\leq \frac12$ for all $x\geq M+1$),
	\begin{align}
		\bE e^{ \beta\sum_i \sigma_i} &= \prod_{i=1}^a \bE e^{ \beta\sigma_i}  
		= \prod_{i=1}^a \bE \frac{\xi_{X'_i}^{-1}}{\xi_{X'_i}^{-1}-\beta} 
		= \prod_{i=1}^a \bE \frac{1}{1-\beta\xi_{X'_i}}  \\
		&= \prod_{i=1}^a \frac1N\sum_{x=1}^N \frac{1}{1-\beta\xi_{x\lor (M+1)}} \\
		&\leq \prod_{i=1}^a \frac1N\sum_{x=1}^N (1+2\beta\xi_{x\lor (M+1)}) 
		= \prod_{i=1}^a \big(1+\frac{2\beta}{N}\sum_{x=1}^N \xi_{x\lor (M+1)}\big)    
	\end{align}
	By \eqref{eq:xilambdasum} and \eqref{eq:xipprox}, 
	\begin{align}
		\sum_{x=1}^N \xi_{x\lor (M+1)} \leq M^{1-\frac1\alpha} + C(\xi)M^{1-\frac1\alpha} = (1+C(\xi))M^{1-\frac1\alpha}.
	\end{align}
	Hence, 
	\begin{align}\label{eq:expo}
		\bE e^{ \beta\sum_i \sigma_i} 
		&\leq \prod_{i=1}^a (1+2\beta(1+C(\xi))M^{1-\frac1\alpha}N^{-1})    
		\leq e^{2\beta(1+C(\xi)) a M^{1-\frac1\alpha} N^{-1} }.
	\end{align}
	Assuming $t\geq 2(1+C(\xi))a M^{1-\frac1\alpha} N^{-1}$, we pick the maximal possible value for $\beta$, namely $\beta = \frac12\xi_{M+1}^{-1}$. Plugging  \eqref{eq:expo} into \eqref{eq:cheby}, we get
	\begin{align}
		\bP(J'_t < a ) \leq e^{-\frac12\xi_{M+1}^{-1}( t-  2(1+C(\xi)) a M^{1-\frac1\alpha} N^{-1} )}.
	\end{align}
	Choosing $a = \frac{N M^{\frac1\alpha-1} }{4(1+C(\xi))} t$ it follows that 
	\begin{align}
		\bP\Big(J'_t < \frac{N M^{\frac1\alpha-1}}{4(1+C(\xi))} t \Big) \leq e^{-\frac14\xi_{M+1}^{-1}t}
		\leq e^{-\frac14C(\xi)M^{\frac1\alpha}t}.
	\end{align}
	for all $t\geq 0$ and 
	\begin{align}
		\int_0^\infty \bP\Big(J'_t < \frac{N M^{\frac1\alpha-1} }{4(1+C(\xi))} t \Big) \;\di t 
		\leq \frac4{C(\xi)}M^{-\frac1\alpha}.
	\end{align}
\end{proof}

\subsection*{The proof of Theorem~\ref{thm:comingdown}}
We conclude the section with the final derivation of the `coming down from infinity'-property.
\begin{proof}[Proof of Theorem~\ref{thm:comingdown}]
By Proposition~\ref{thm:coalescence-time} and Fatou's lemma,
\begin{equation}\label{eq:limas}
	\bE \tau_{\infty\to1} = \sum_{M=2}^\infty \bE \tau_{M+1\to M} 
	\leq \sum_{M=2}^\infty \liminf_{N\to\infty} \bE \tau^\super{N}_{M+1\to M}
	\leq C(\xi)\sum_{M=2}^\infty M^{-\frac1\alpha}
	<\infty.
\end{equation}
In particular, $\tau_{\infty\to1}$ is a.s. finite. Let $\mathcal{N}(t)=\sum_{x\in\mathbb{N}} \gamma_t(x)$ the process counting the number of particles at time $t$. By \eqref{eq:limas} infinity is a continuous entrance point. By Lemma 1.2.\ together with Remark 2.3.\ of \cite{Foucart20}, the finiteness
\(
\sup_{m}\bE[\tau_{m\to1}]
\)
implies that $(\mathcal{N}(t))_{t\geq 0}$ can be extended to a Feller process over $\bN\cup\{\infty\}$ such that under the extended law, $\mathcal{N}$ starts from $\infty$-many particles and instantaneously enters a configuration with only finitely many particles.  The analogous statement for the whole process $\gamma$ now follows from the Feller property of the particle count process together with Proposition~\ref{prop:projlim} and Lemma~\ref{lem:Feller} as well as the fact that we may encode the particles that reside at $\infty$ at any time as an element of $2^{\bN}$.
\end{proof}

\section{Proofs of the main results}
In this section, we finally prove Theorem~\ref{thm:limitvotermodel}, Theorem~\ref{thm:limitvotermodel2} and Theorem~\ref{thm:votermodelconsensus}.
\begin{proof}[Proof of Theorem~\ref{thm:limitvotermodel}]
In Proposition~\ref{prop:convYtoX}, it is shown that the processes $Y^{(N)}$ under temporal rescaling by $a_N$ converge in distribution to a limit random walk $X$ on $\bN$ which is governed by rates $\xi$ as in \eqref{eq:limweights}. Proposition~\ref{prop:convYtoXfinite} extends this to systems $(\gamma_t)_{t\geq 0}$ of finitely many coalescing random walks. Hence, for any $t\geq 0$ and any $A\subset \bN$ finite we may define the corresponding probabilities for the limiting voter model through
\begin{equation}\label{eq:duality2}
	\bP\big(\eta_t(x)\equiv 1 \text{ on }A\big| \eta_0\big):=\bP\big(\eta_0(x)\equiv 1 \text{ on }\supp \gamma_t\big| \supp \gamma_0 = A\big),
\end{equation}
where $\supp \gamma_t=\{x:\gamma_t(x)=1\}, t\geq 0.$ The extension of the definition to arbitrary $A\subset \bN$ can now be achieved by observing that the dual limiting system is well defined for an arbitrary initial configuration of particles dispersed over $\bN$ by Proposition~\ref{prop:projlim}.
\end{proof}
\begin{proof}[Proofs of Theorems~\ref{thm:limitvotermodel2} and~\ref{thm:votermodelconsensus}]
We only derive Theorem~\ref{thm:limitvotermodel2} explicitly; the statements for the finite system can be obtained in the same fashion but without relying on properties of the limiting random walk $X$. By duality, the expectation of the consensus time $\tau$ is clearly bounded above by the expectation of the time $\tau_{\infty\rightarrow 1}$ that is needed in the dual system to coalesce the initial configuration with a particle at every site into one single particle. By Theorem~\ref{thm:comingdown}, this time has a finite mean, which proves the first part of Theorem~\ref{thm:limitvotermodel2}. To prove formula~\eqref{eq:fixationprobformula}, note that at any fixed time $t_0>0$, we may pick a site $x$ with probability proportional to $\xi_x$ and obtain the type $\eta_{t_0}(x)=\eta_0(y)$ by tracing the path of the corresponding particle in the dual for time $t_0$ to some target site $y$. The distribution of this path is independent of whether consensus has occurred and agrees with the distribution of a path of the limiting walk $X$ initiated from its invariant distribution by Lemma~\ref{lem:stationary}. In particular, $y$ has the invariant distribution $(\xi_x/\sum_z\xi_z)_{x\in\bN}$ of $X$ if $t_0\geq \tau$, which leads immediately to~\eqref{eq:fixationprobformula}.
\end{proof}

\section*{Acknowledgments}
We thank Florian V\"ollering for early collaboration and contributions to this project.

\section*{Funding}
This collaborative research was initiated within the scientific network \emph{Stochastic processes on evolving networks} funded by Deutsche Forschungsgemeinschaft (DFG, German Research Foundation) -- grant no. 412848929. CM's research was funded by Deutsche Forschungsgemeinschaft (DFG, German Research Foundation) -- grant no. 443916008 (SPP 2265).

\bibliographystyle{plainnat}
\bibliography{AVM}

\end{document}